\documentclass[12pt]{article}

\usepackage{amssymb}
\usepackage{graphicx}
\usepackage{amscd,amssymb,amsfonts}
\usepackage{color}
\usepackage{eurosym}
\usepackage[english]{babel}
\usepackage[utf8]{inputenc}
\usepackage[T1]{fontenc}
\usepackage{epsf}
\usepackage{amsmath}
\usepackage{amssymb}
\usepackage{graphics}
\usepackage{mathtools}
\usepackage{array}
\usepackage{latexsym}
\usepackage{graphicx}
\usepackage{color}
\usepackage[table]{xcolor}
\newcommand{\R}{{\mathbb R}}
\usepackage{url}

\begin{document}
\pdfoutput=1
\title{Shape recovery from sparse tomographic X-ray data}

\author{Heikki Haario, Aki Kallonen, Marko Laine, Esa Niemi, Zenith Purisha \and Samuli Siltanen}

\maketitle

\begin{abstract}
A two-dimensional tomographic problem is studied. The target is assumed to be a homogeneous object bounded by a smooth curve. A Non Uniform Rational Basis Splines (NURBS) curve is used as computational representation of the boundary. This approach conveniently provides the result in a format readily compatible with computer-aided design (CAD) software. However, the linear tomography task becomes a nonlinear inverse problem due to the NURBS-based parameterization. Therefore, Bayesian inversion with Markov chain Monte Carlo (MCMC) sampling is used for calculating an estimate of the NURBS control points. The reconstruction method is tested with both simulated data and measured X-ray projection data. The proposed method recovers the shape and the attenuation coefficient significantly better than the baseline algorithm (optimally thresholded total variation regularization), but at the cost of heavier computation.
\end{abstract}

\newpage
\tableofcontents

\newpage

\section{Introduction}\label{Introduction}

We propose a new reconstruction algorithm for two-dimensional X-ray tomography of objects with homogeneous attenuation coefficient. The basic idea is to represent the boundary of the object in terms of a non-uniform rational basis spline (NURBS) curve, thereby providing a direct connection to computer-aided design (CAD) software \cite{Bazilevs,Cottrell,Piegl,Rogers2001}. The control points of the parametric curve are the degrees of freedom in the inverse problem, which then becomes nonlinear. Therefore, we resort to Bayesian inversion and Markov chain Monte Carlo (MCMC) sampling \cite{Kaipio,Kolehmainen,Robert,Siltanen} for estimating the unknown curve from projection data.

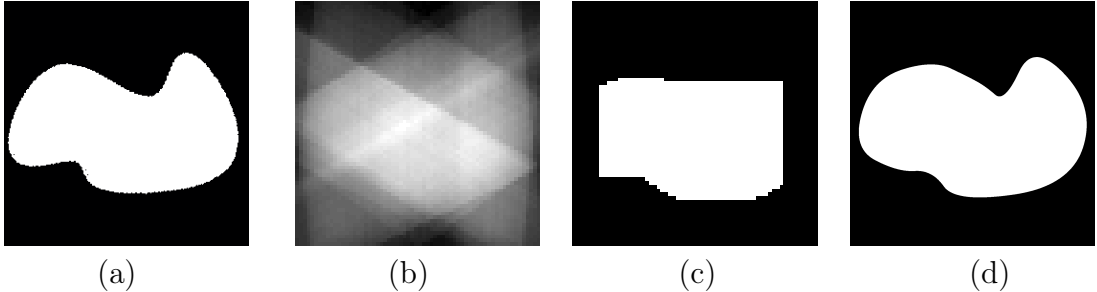
\begin{figure}[t!]
\begin{picture}(200,120)
\put(-15,4){\includegraphics[width=3.25cm]{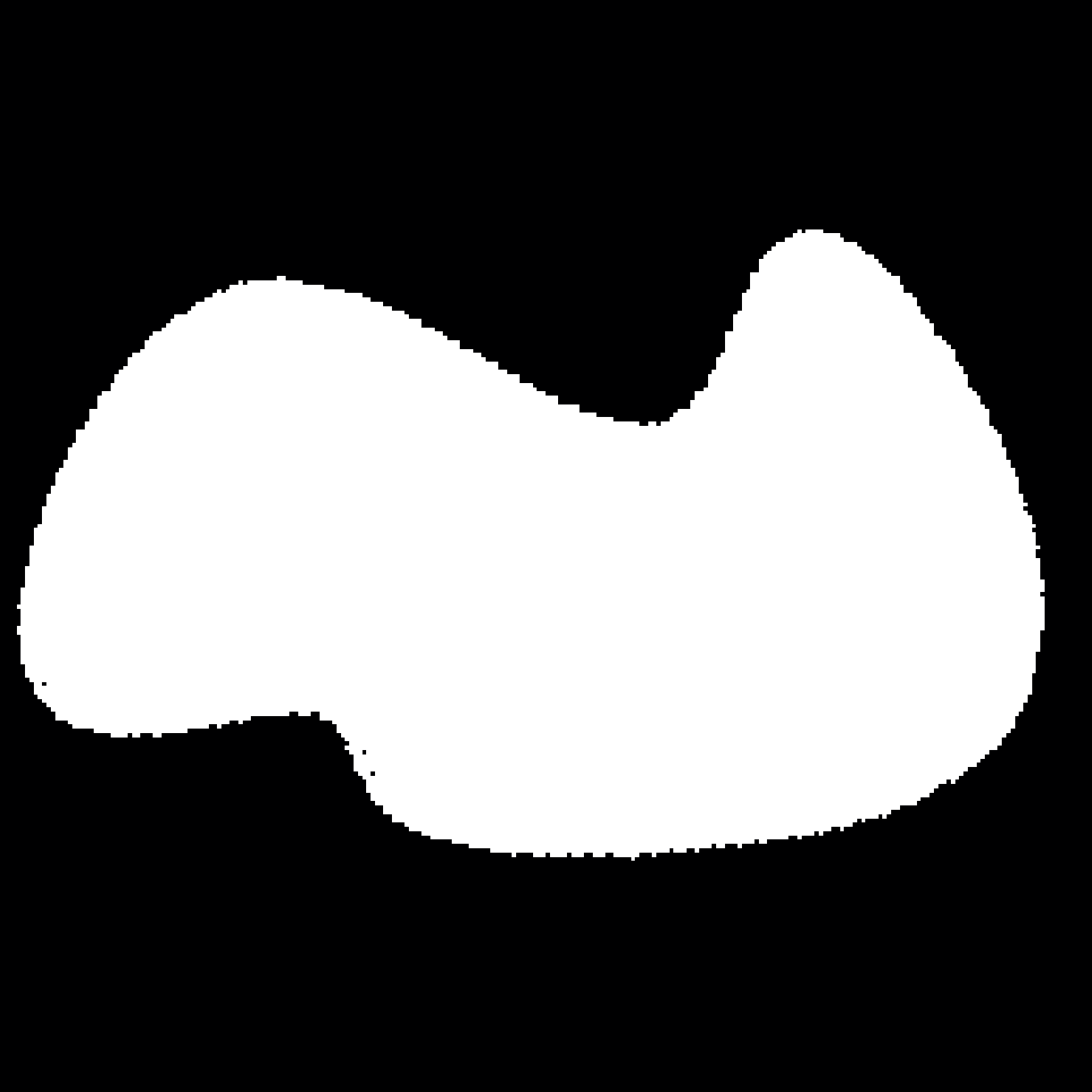}}
\put(95,4){\includegraphics[width=3.25cm]{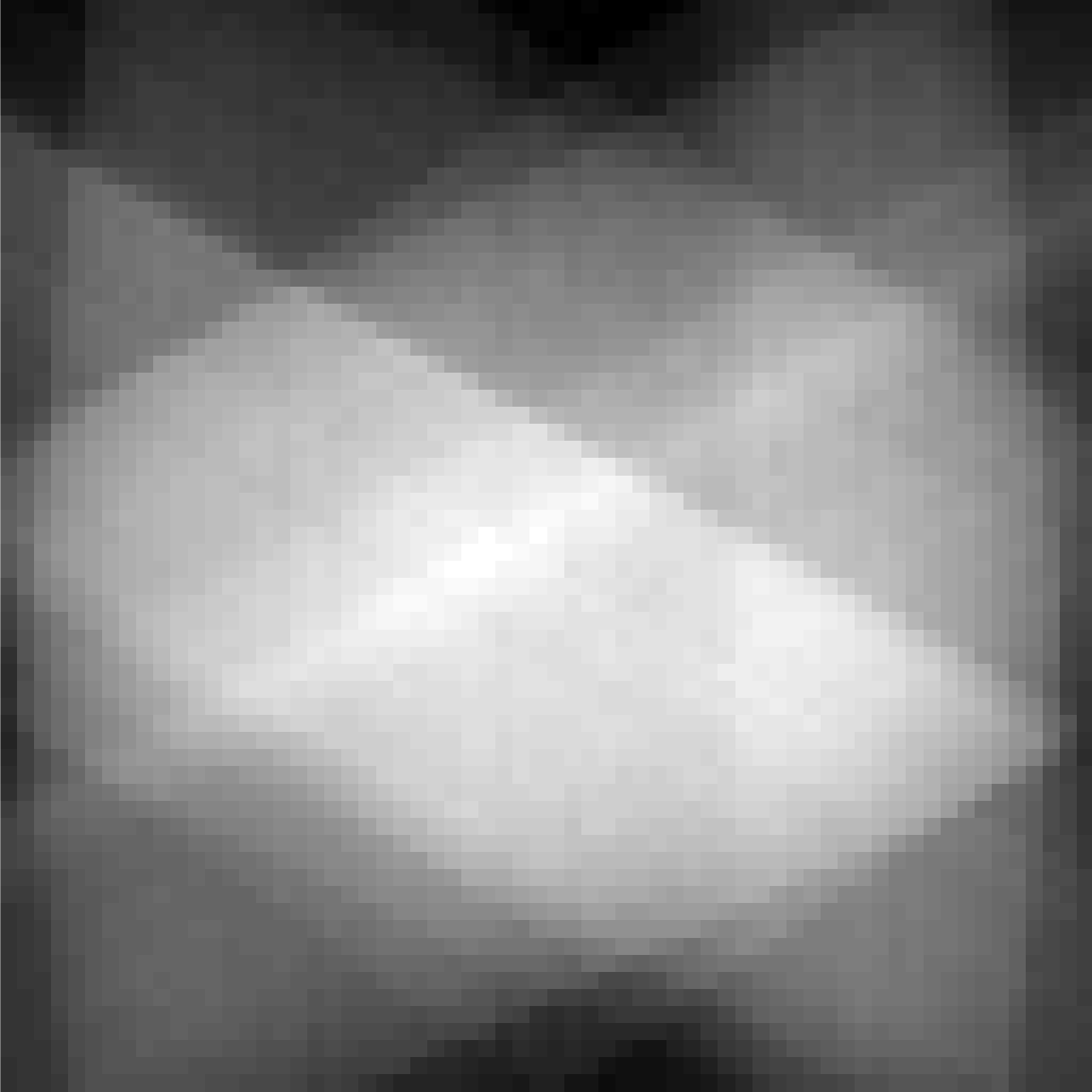}}
\put(200,4){\includegraphics[width=3.25cm]{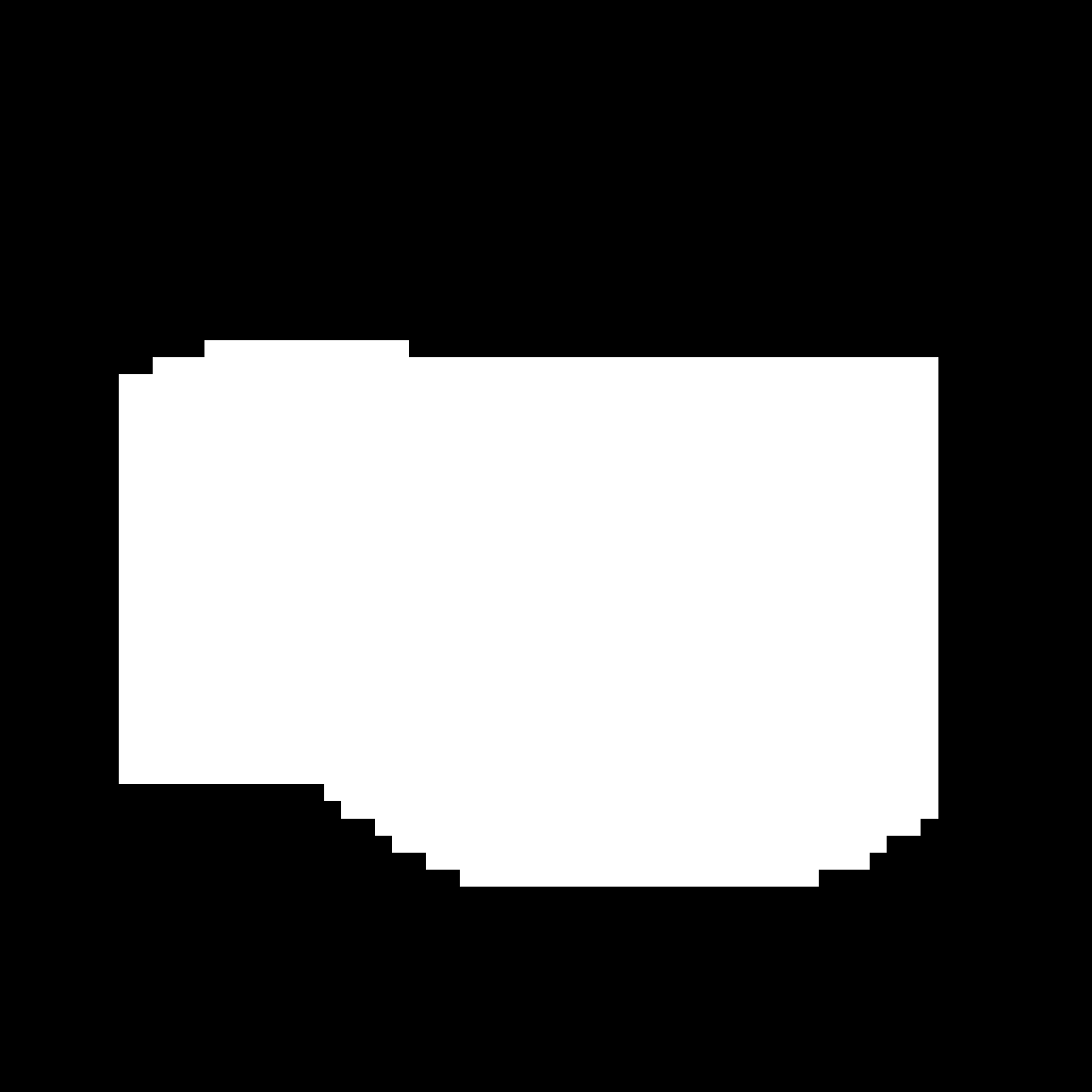}}
\put(305,4){\includegraphics[width=3.25cm]{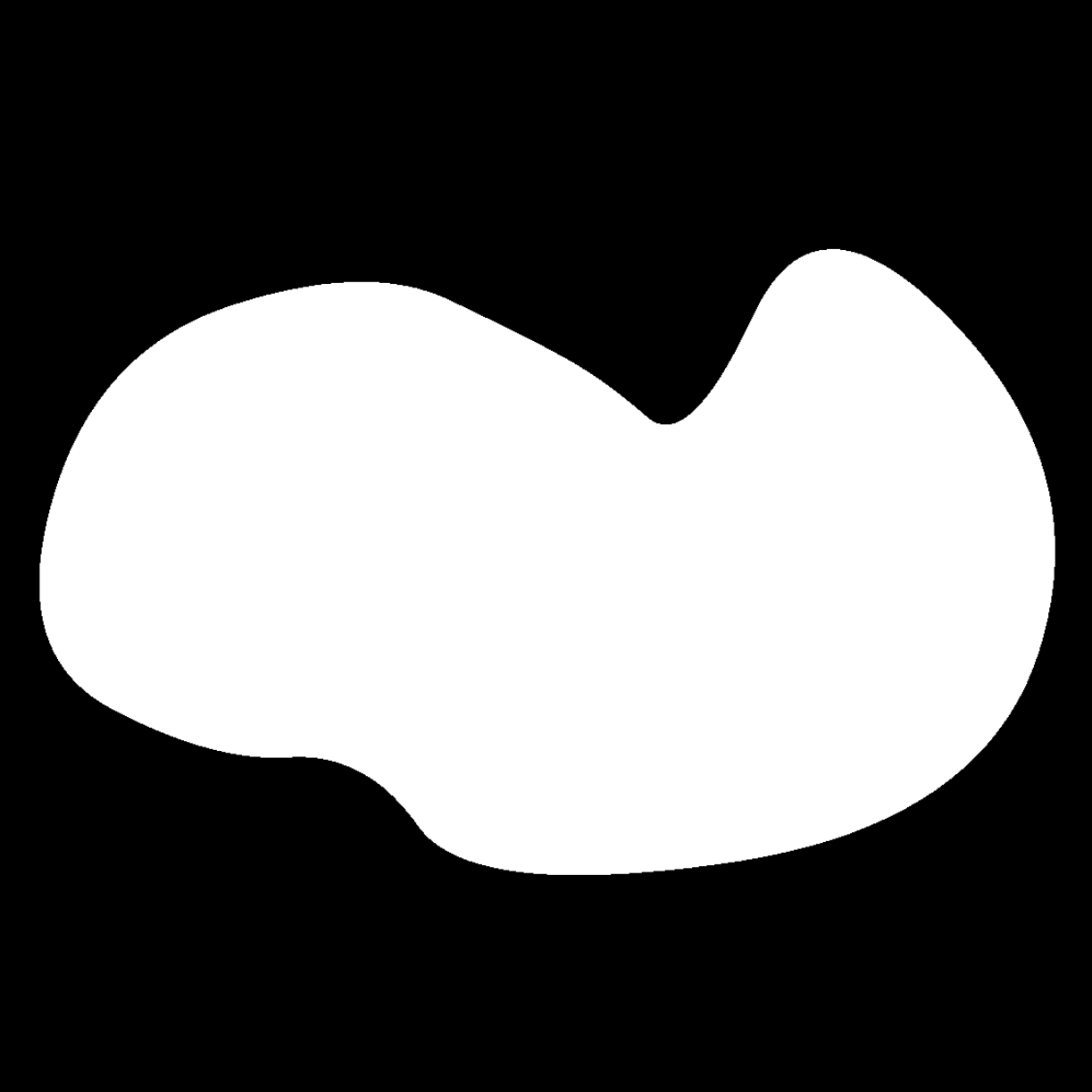}}
\put(20,-10){(a)}
\put(130,-10){(b)}
\put(240,-10){(c)}
\put(350,-10){(d)}
\end{picture}
\caption{\label{fig:hammer} (a) shape of the true object. (b) filtered back projection reconstruction from 6 measured X-ray projections. (c) optimally thresholded total variation (TV) regularized reconstruction from the same data. The regularization parameter was optimally chosen. The result is given as a pixel image. (d) MCMC-NURBS estimate computed from the same data. The result is given in CAD-compatible vector-graphic format.}
\end{figure}

Why not reconstruct the object with a traditional reconstruction method, such as filtered back-projection, and fit a NURBS curve to the boundary of the reconstructed domain? This indeed is quicker and more reliable in the case of  a comprehensive tomographic dataset with dense angular sampling. However, often there are time constraints, geometric obstructions, or radiation dose issues preventing the collection of  a detailed dataset. With such cases in mind, we test our algorithm with projection data collected from only {\em six} directions around the object. It turns out that the combination of NURBS and MCMC outperforms many traditional methods in this extremely ill-posed inverse problem. See Figure \ref{fig:hammer}.

The proposed MCMC-NURBS reconstruction method is first tested with two simulated phantoms, one of them convex and the other non-convex. See \cite{Purisha2014,Purisha2016} for other examples. The X-ray attenuation coefficient is constant inside the phantoms, and we simulate six projection images using fan-beam geometry. The projection directions are uniformly distributed over the full circle. We use Delayed Rejection Adaptive MCMC (DRAM) approach for effective Monte Carlo computation \cite{Andrieu,Haario,Robert}. The degrees of freedom in the inverse problem are the coordinates of the control points of the NURBS curve and the unknown attenuation coefficient.

We use Total Variation (TV) regularization as a baseline method for comparison and quantitative assessment of tomographic reconstruction quality. The degrees of freedom in the inverse problem are the non-negative pixel values. To be as fair as possible to the baseline method, we choose the regularization parameter optimally by comparing the reconstruction to the true simulated phantom. We force the reconstruction to be piecewise constant (with only one possible positive attenuation value) by optimally thresholding the TV regularized solution. 

Furthermore, we test both methods using X-ray projection data measured from physical objects, similarly shaped than the simulated phantoms. 

The proposed method is found to outperform the optimized TV regularization algorithm in all the cases studied. See Figure \ref{fig:hammer} for the example concerning the physical non-convex phantom.

NURBS curves have been used as building blocks in  CAD software for decades. Some works in optimizing the NURBS representation in certain inverse problems, including reverse engineering, fitting strategies, and recovering shapes from photographs, have been done in \cite{Alpers,Anson,Jing,Khan,Ma,Renken,Sarfraz,Saini}.  Other works in X-ray tomography that readily provide segmented images can be found in \cite{Djafari2,Djafari1,Klann}. Non-destructive testing from very restricted  data using Bayesian inversion in terms of a vector-graphic format are discussed in \cite{Djafari2,Djafari1,Kolehmainen}.

Also, there are various branches of computational science using NURBS curves representation.  For example, isogeometric analysis  provides a computational approach for integrating finite element analysis (FEA) and CAD \cite{Cottrell}. One of the applications of NURBS-based isogeometric analysis is computing flows from spinning propellers in opposite direction. In this case, a discretization using NURBS gives more accurate results in computing the flows compared to standard finite elements \cite{Bazilevs}. 

The method proposed here is different from all the previous ones listed above.

\section{Materials and Methods}\label{Materials and Methods}
\subsection{NURBS}\label{NURBS}

We consider an unknown physical domain $ \Omega \subset  \R^2$ and model the object boundary $\partial \Omega$ defined by a parametrized curve $\mathcal{S}:[0,1]\rightarrow \mathcal{S} ([0,1]) = \partial \Omega \subset\R^2$. In our computational problem, a NURBS curve as a piecewise rational function $\mathcal{S}$ defined on $t\in [0,1]$, is introduced. We divide the interval [0,1] into $K-1$ pieces $0 = t_1 \le t_2 \le ... \le t_{K} = 1$ and called {\it breakpoints} which map into the endpoints of the polynomial segment.

We introduce the curve representation as follows:

\begin{equation}
\mathcal{S}(t)=\sum_{i=0}^n\textbf{p}_iR_{i,p}(t), \quad 0\leq t \leq 1 \\
\end{equation}
where the $n+1$ points $\textbf{p}_i\in\R^2$ configure the curve shape. The $n+1$ points are called control points and $R_{i,p}(t)$ is the following rational function of degree {\it p}:

$$
R_{i,p}(t)=\frac{\omega_iN_{i,p}(t)}{\sum_{i=0}^n\omega_iN_{i,p}(t)},
$$
where $\{\omega_i\}$ are nonnegative weights for all $i$ and  $\{N_{i,p}(t)\}$ are basis functions that describe how strongly the control points $\{\textbf p_i\}$ attract the NURBS curve. They are defined recursively as

\begin{equation*}
N_{i,0}(t)=
\begin{cases}
1 & \text{if $t_i\le t  < t_{i+1},$}\\
0 & \text{otherwise,}
\end{cases}
\end{equation*}
and for $p>0$ as
\begin{eqnarray}
N_{i,p}(t)=
\frac{t-t_i}{t_{i+p}-t_i}N_{i,p-1}(t)+ 
\frac{t_{i+p+1}-t}{t_{i+p+1}-t_{i+1}}N_{1+i,p-1}(t),
\end{eqnarray}
where $\frac{0}{0}=0$ by definition.
A collection of $K$ {\it breakpoints} is then called knot vector:

\begin{equation}
{\mathbf t} = [ t_1, t_2, ..., t_{K} ]^T,
\end{equation}
where $K = p + n + 2$.

If $0 = t_1 < t_2 < ... < t_{K} = 1$ are evenly spaced, then ${\mathbf t}$ is called a periodic uniform knot vector. 
In our discussion, we implement a closed NURBS curve using periodic uniform knot vector to recover the boundary of the object. By repeating the first $p$ control points after the last point, an unclamped closed NURBS curve is obtained. 

\begin{figure}[h!]
\begin{picture}(125,125)
\put(-10,-15){\includegraphics[width=7cm]{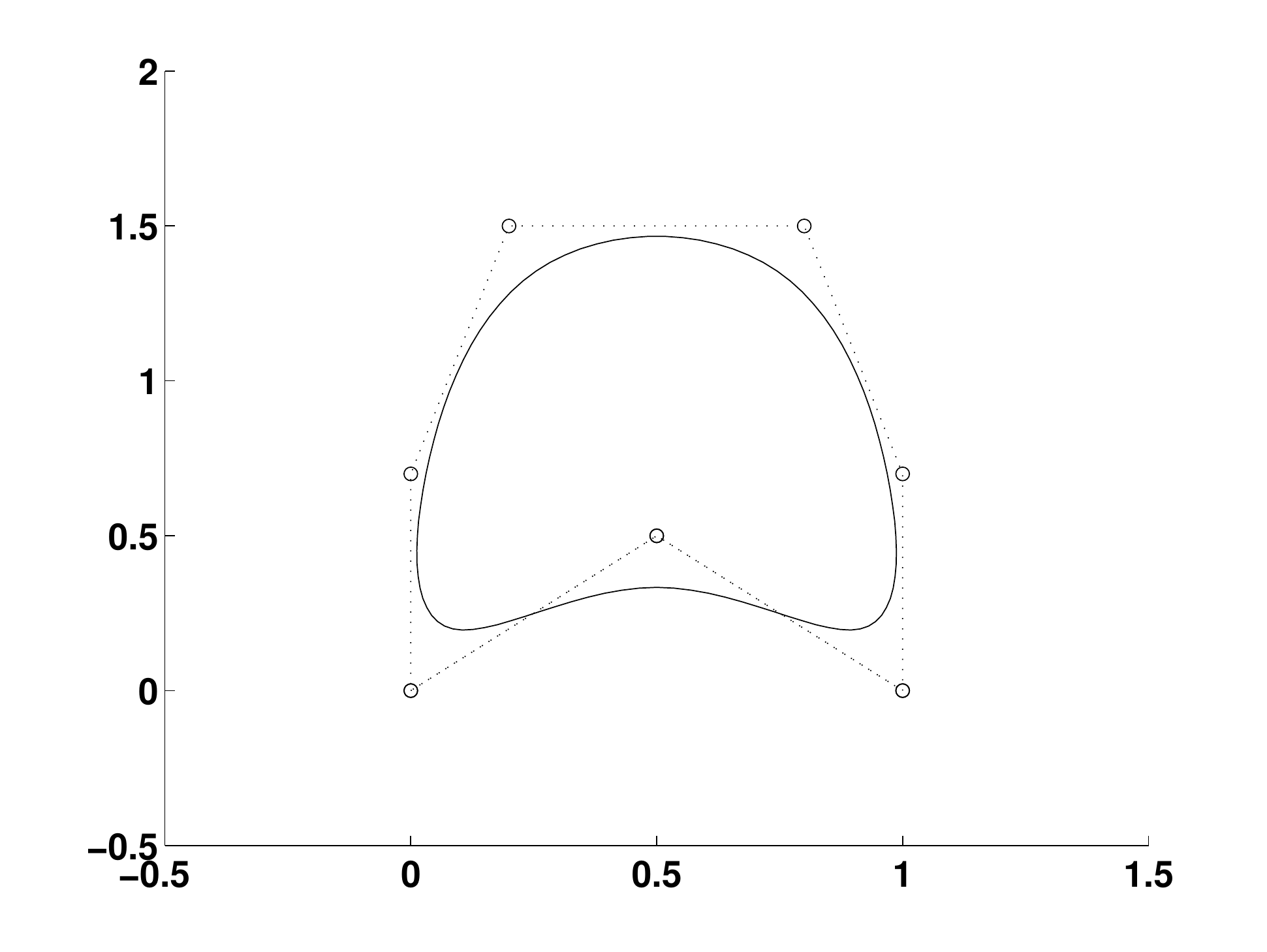}}
\put(180,-15){\includegraphics[width=7cm]{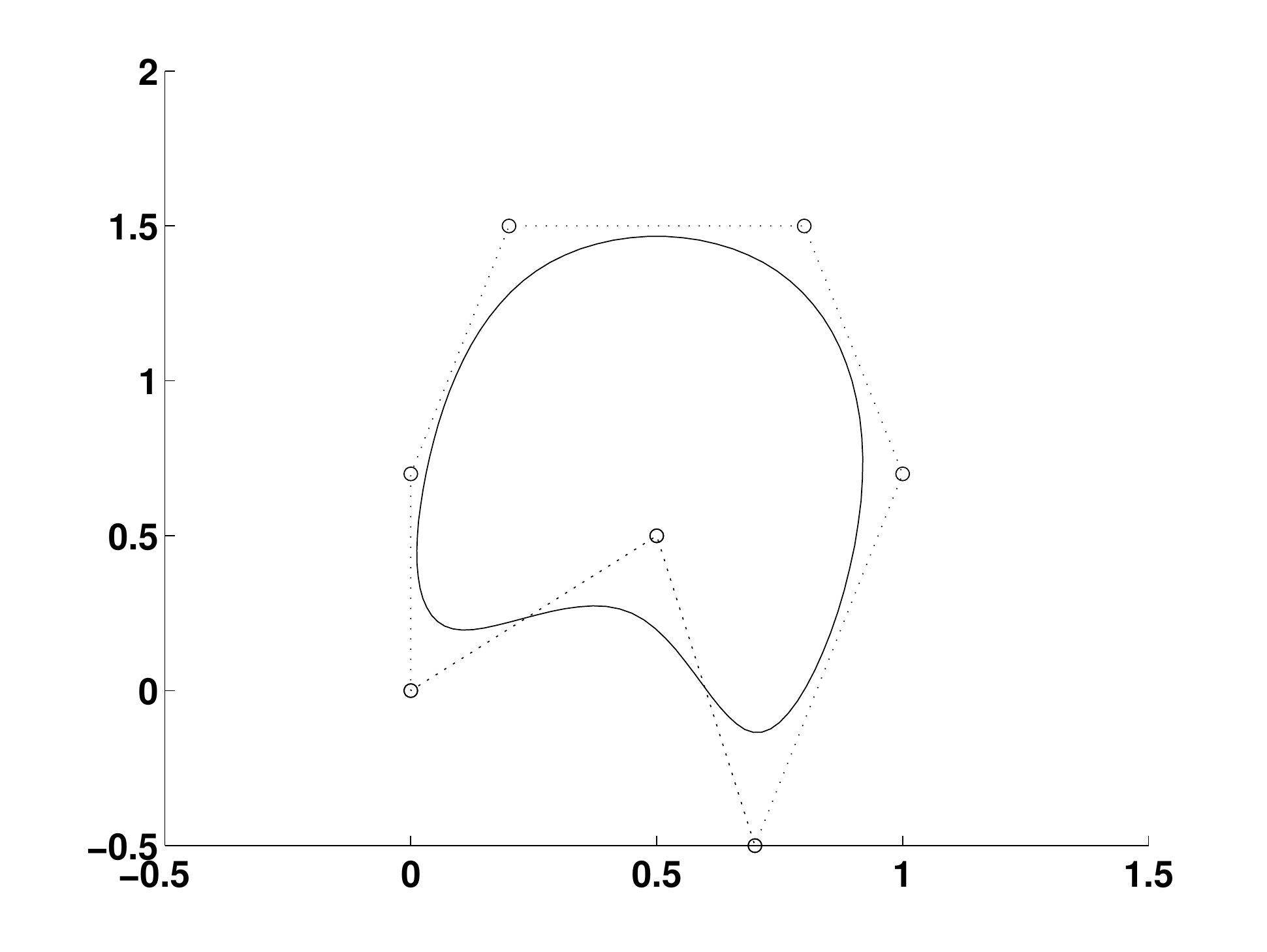}}
\end{picture}
\caption{\label{fig:closedcurve}Left: NURBS curve with control points $(0,0),$ $(0.5,0.5),$ $(1,0),$ $(1,0.7),$ $(0.8,1.5),$ $(0.2,1.5),$ $(0,0.7)$. Right: NURBS curve with control points $(0,0),$ $(0.5,0.5),$ $(0.7,- 0.5),$ $(1,0.7),$ $(0.8,1.5),$ $(0.2,1.5),$ $(0,0.7)$.}
\end{figure}

Two examples of closed NURBS curves of $3$rd-degree basis functions with seven control points are given. Figure \ref{fig:closedcurve} illustrates the closed curve with the basis functions defined on the same periodic uniform knot vector:
$$
[0, \frac{1}{13},  \frac{2}{13},  \frac{3}{13},  \frac{4}{13},  \frac{5}{13},  \frac{6}{13},  \frac{7}{13},  \frac{8}{13},  \frac{9}{13},  \frac{10}{13},  \frac{11}{13},  \frac{12}{13}, 1]^T.
$$

In those examples and in our computations, we assume that the weights corresponding to all of the control points are the same.

\subsection{Tomographic Measurement Model}\label{Tomographic Measurement Model}

Consider a physical domain $\Omega \subset \R^2$ and a continuous tomography model $f:\R^2\rightarrow\R$ where $f(x,y) \ge 0$ as 

\begin{equation}\label{eq:f}
 f(x,y) =
\begin{cases}
c,    \mbox{ for } (x,y)\in\Omega\\
0,   \mbox{ for } (x,y)\in\R^2\setminus\Omega,
\\
\end{cases}
\end{equation}
and $supp(f) \subset \Omega$.

\begin{figure}
\begin{picture}(100,100)
\put(90,-15){\includegraphics[width=5cm]{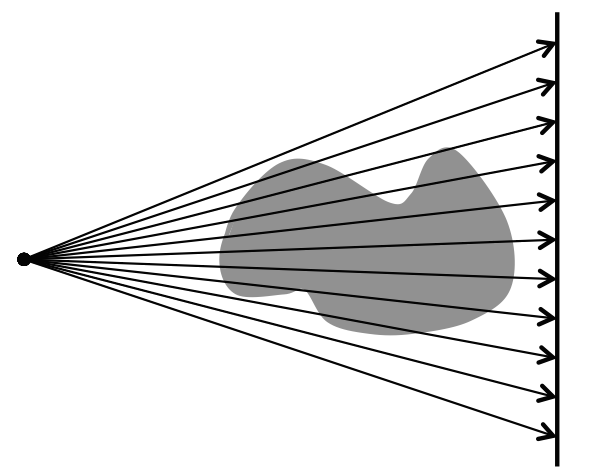}}
\end{picture}
\caption{\label{fig:fan beam}Fan beam X-ray measurement geometry. The dot represents the locations of the X-ray source. The arrows show the detector measuring the intensity of the X-rays after passing through the target.} 
\end{figure}

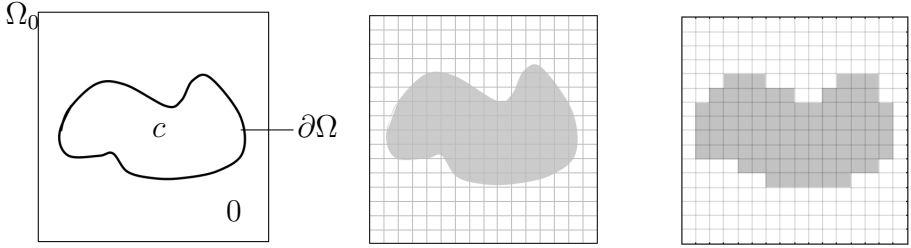
\begin{figure}
\begin{picture}(90,90)
\put(15,7){\includegraphics[width=3.25cm]{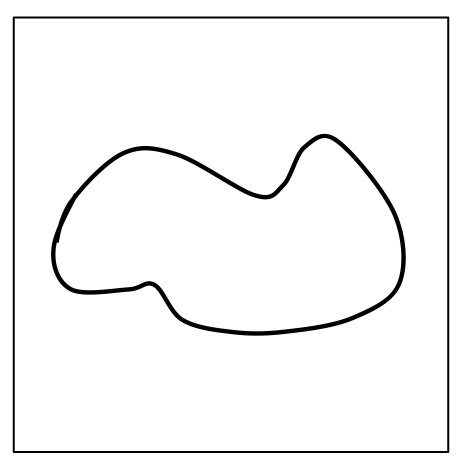}}
\put(94,52){\line(1,0){20}}
\put(140,6){\includegraphics[width=3.25cm]{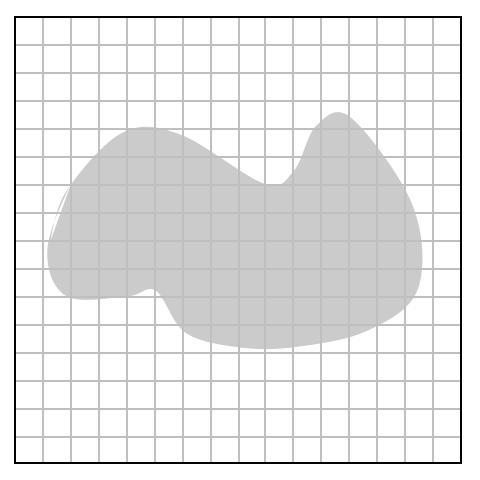}}
\put(5,92){$\Omega_0$}
\put(115.5,49.5){$\partial\Omega$}
\put(60.5,49.5){$c$}
\put(88.5,17.5){$0$}
\put(260,8){\includegraphics[width=3.1cm]{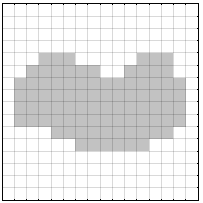}}
\end{picture}
\caption{\label{fig:finite models} Left: the NURBS-based object model where $\partial\Omega$ is a NURBS curve. Middle: The attenuation function $f$ defined in (\ref{eq:f}). The ${\mbox{N $\times$ N}}$ pixel grid is indicated in gray. Right: the pixel-based attenuation model (\ref{eq:B}).}
\end{figure}

For computational reasons, we need to construct two discrete models: a pixel-based object model and a NURBS-based object model as in Figure \ref{fig:finite models}. Pixel-based object model is a discretization of the continuous model into pixels in two dimensional image, so then it discretizes a line integral to a standard pencil-beam model. NURBS-based object model is a model where the boundary of the object $\partial \Omega$ is represented in terms of a NURBS curve. 
We represent the control points $(x_i,y_i), i=0,...,n$ in polar coordinates where ${\bf n}$ is the number of control points:
 $$
 x_i=r_i\cos\theta_i ; y_i=r_i\sin\theta_i.
 $$
 
Consider a vector $v\in\R^{2n+3}$ and construct $v=[r_0, \theta_0, ... ,r_n,\theta_n,c ]^T$, where $v_{2n+3}=c$ is an attenuation parameter.

The line integral is discretized by switching to pixel-based model using an operator 
  $\mathcal{B}:\R^{2n+3} \rightarrow \R^{\mbox{\tiny {N$\times $N}}}$ where $\mbox{N$\times$N}$ is the resolution of the pixel image. Define

\begin{equation} \label{eq:B}
\mathcal{B}(v)_{ij}  =
\begin{cases}
c,   \text{ if the pixel center $(i,j)$ is inside the NURBS curve,}\\
0,  \text{ if the pixel center is outside the NURBS curve.}
\end{cases}
\end{equation}

\begin{figure}[h!]
\begin{picture}(100,100)
\put(120,0){\includegraphics[width=3.25cm]{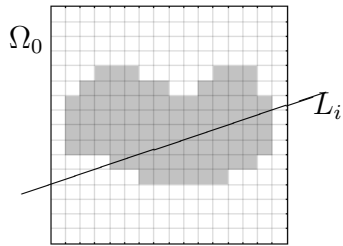}}
\put(110,20){{\line(3,1){115}}}
\put(105,75){$\Omega_0$}
\put(220,50){$L_i$}
\end{picture}
\caption{\label{fig:line integral}The standard pencil-beam model is applied to the pixel-based model.}
\end{figure}

Let consider $\mathcal{B}(v) \in \R^{\mbox{\tiny {N$\times $N}}}$ as in Figure \ref{fig:line integral}.
In the pixel-based model, the measurement is defined as follows:
$$ 
{\it m_i} = \int_{L_i} f(x,y) ds \approx \sum_{j=1}^{n} a_{ij} b_{ij},
$$
where $a_{ij}$ is a distance that the line $L_i$ travels in the j{\it th} pixel. We use the collection of lines $L_i$ arising from several fan-beam projections as shown in Figure \ref{fig:fan beam}. Our measurement model can be written as follows.

\begin{equation} \label{eq:m}
{\bf m}=\mathcal{A}(\mathcal{B}(v) ). 
\end{equation}

\subsection{Tomographic Projection Data}\label{Tomographic Projection Data}
The x-ray tomography data of the sugar phantom were acquired with the custom-built $\mu$CT device nanotom 180 supplied by Phoenix| \ Xray Systems + Services GmbH (Wunstorf, Germany).  The chosen geometry resulted in a magnification with resolution of 62.2 $\mu$m/pixel for physical phantom $\Omega_1$ and 63.8 $\mu$m/pixel for physical phantom $\Omega_2$. The phantoms can be seen in Figure \ref{fig:Physical phantoms}. The x-ray detector is a 12-bit CMOS flat panel detector with $1128 \times 1152$ pixels of 100 $\mu$m size (Hamamatsu Photonics, Japan). A set of 120 projection images were acquired over a full 360 degree rotation with uniform angular step of 3 degrees between projections. Each projection image was composed of an average of six 250 ms exposures. The x-ray tube acceleration voltage was 80 kV and tube current 300 $\mu$A, and the full polychromatic beam was used for image acquisition. For this work we choose the projections corresponding to the middle cross-section of a sugar phantom as Figure \ref{fig:Physical phantoms}, thus the task is only 2D problem. We picked 6 projections from the measured data with uniform angular sampling from a total opening angle of 360 degrees.

\subsection{Bayesian Inversion}\label{Bayesian Inversion}

In our problem, we have a sparse tomographic model which leads to ill-posed inverse problem. To compensate this ill-posedness, {\it a priori} knowledge of the target should be explored well. For computational reasons, the information needs to be transformed to {\it quantitative} form. Bayesian inversion provides a flexible way to synthesize this extra information of the target \cite{Kaipio}.

The main idea of this approach is to recast the inverse problem as a problem of Bayesian inference. We use probability theory to model our lack of information in the inverse problem. All the variables in the model are considered as random variables.

There are three main steps in constructing the problem to Bayesian framework. Firstly, gathering information prior to the measurement and construct it as a prior density. Secondly, constructing a likelihood function that expresses how likely the observation outcomes with unknown parameter given. Thirdly, constructing and exploring the posterior probability density as a solution of the inverse problem.	

In our case, the tomographic model leads to the nonlinear inverse problem with additive Gaussian errors $\varepsilon$  as follows:
 
\begin{equation}\label{eq:BI2}
{\bf m}=\mathcal{A}(\mathcal{B}(v) )+\varepsilon. 
\end{equation}
where $\mathcal {A}$ is a matrix model of the measurement and \textbf{m} is the x-ray measurement.

\subsubsection{Prior distribution} 

\begin{figure}
\begin{picture}(100,120)
\put(100,-15){\includegraphics[width=5cm]{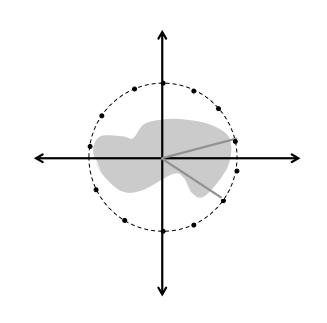}}
\put(185,65){r}
\put(184,50){$\theta$}
\end{picture}
\caption{\label{fig:prior}The position of the control points in the prior is shown.}
\end{figure}

In this part, we construct a prior density which will contribute to remove the ill-posedness. In this case, we have a homogeneous object with diameter not more than $2r$, where $r$ is determined as a radius from the origin. This condition means that the control points satisfy $|p_i|\leq2r$. They are shown as the dots in the Figure \ref{fig:prior}. The dashed circle is the NURBS curve produced by the control points with a fixed uniform periodic knot vector.

Assume that the prior information is Gaussian distributed with variance $\sigma_0$. Then a {\em priori} information has the following quantitative form

 \begin{equation}
\pi(v)= \exp(-\frac{1}{2\sigma_{0}^2}\|v -\widetilde v \| _2^2),
\end{equation}
 where
$v=[r_1, \theta_1, ... ,r_n,\theta_n,c ]^T$, $ \widetilde {v} =[\widetilde{r_1}, \widetilde{\theta_1}, ... ,\widetilde{r_n},\widetilde{\theta_n},\widetilde{c} ]^T$.

\begin{figure}
\begin{picture}(110,100)
\put(150,-10){\includegraphics[width=2.5cm]{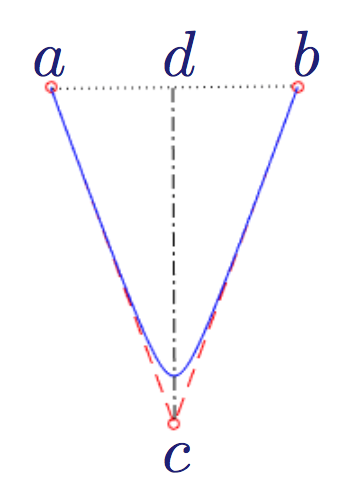}}
\end{picture}
\caption{\label{osc prior}The illustration of hard prior is given, where $a,b,{\text {and }} c$ are the control points.}
\end{figure}

In addition, to avoid uncontrollable movement of control points, where they can rotate only to one side and get trapped, we condition the range of the control points movement by setting $\theta_i$ as an angle of each control point: 
\begin{equation}\label{eq:boundary_prior}
\max \big \{ \theta_{i-2}, \Gamma_i^{\text min} \big \} \leq \theta_i \leq \min \big \{ \theta_{i+2}, \Gamma_i^{\text max} \big \},
\end{equation}
where $\Gamma_i^{\text min}$ is a lower bound for $\theta_i$ and $\Gamma_i^{\text max}$ is an upper bound for $\theta_i$.  The condition of the radius length is provided as well as: $0 \leq r_i \leq r^{\mbox {\tiny M}}$ where $r^{\mbox {\tiny M}}$ is a constant.
To minimize oscillations in the curve, the following condition should be satisfied

$$
||d-c|| \leq k ||a-b||
$$ 
where $a,b,{\text {and }} c$ are the control points as in Figure \ref{osc prior} and $k$ is a constant. 

To avoid a control polygon intersection, another hard prior setting is added. Assume that $ef$ is another polygon in the Figure \ref{osc prior}, then the self-intersection can be avoided if $(\vec {ea} \times \vec {ef}).(\vec {ec} \times \vec {ef}) \leq 0$ holds.

\subsubsection{Likelihood function} 
We have $\varepsilon\sim\mathcal{N}(0,I\sigma^2)$, so then

$$
 {\bf m} - \mathcal{A}(\mathcal{B}(v))  \sim\mathcal{N}(0,I\sigma^2).
$$
From Equation \ref{eq:BI2} we then model the measurement process as:
$$
  \pi({\bf m}\,|\,v) = C\exp(-\frac{1}{2\sigma^2}\|\mathcal{A}(\mathcal{B}(v))- {{\bf m}}\|_2^2),
$$ 
and it is called {\em likelihood} function, where $C$ is a normalization constant.

\subsubsection{Posterior distribution} 
The solution of the inverse problem is the {\it posterior} probability distribution which has the following form:
$$
  {\pi(v\,|\, {\bf m}) =  \frac{\pi(v)\pi({\bf m}\,|\,v)}{\pi({\bf m})}}
$$

or

$$
  {\pi(v\,|\, {\bf m}) \sim  {\pi(v)\pi({\bf m}\,|\,v)}},
$$
where $\sim$ denotes equality up to a normalization constant.                                                                                              

To get an estimation of $v$, we use the conditional mean method that is averaging the values of the generated samples:

\begin{equation}\label{def:CM}
v^{\mbox{\tiny CM}} = \int_{\mathbb{R}^{2{\bf n}+1}}\,v \pi(v\,|\,{\bf m})dv.
\end{equation}

Since we face integration problem, Markov chain Monte Carlo (MCMC) technique is proposed to solve this problem. In the following section, details of MCMC method are discussed. 

\subsection{Markov chain Monte Carlo}\label{MCMC}
Let us consider the Equation~\ref{def:CM}. To handle the integration problem, MCMC approach is applied by generating samples from the posterior distribution. In probability theory, the {\it law of large numbers} ensures that:

$$
\frac{1}{N}\sum_{i=1}^{N} v^{(i)} \approx  \int_{\mathbb{R}^{2{\bf n}+1}}\,v \pi(v\,|\,{\bf m})dv.
$$
Hence, the conditional mean (CM) estimate for $v$ can be written as follows:

\begin{equation}\label{def:CMest}
v^{\mbox{\tiny CM}} \approx \frac{1}{N}\sum_{i=1}^{N} v^{(i)}.
\end{equation}
This CM estimate is the result of the recovered control points and attenuation coefficient.

A well known algorithm in MCMC is Metropolis Hastings \cite{Chib,Gelman,Mira,Tierney1994}. The Markov chains are generated as follow:

\begin{enumerate}
\item  For $i=0$, give an initial state ${v}^{(i)}$.
\item  Set the proposal state $v:= v^{(i)}+\epsilon_k$, where $\epsilon_i \sim \mathcal{N}(0,1)$
\item  If $\pi(v|{\bf m}) \ge \pi(v^{(i)} | {\bf m})$ then set $v^{(i+1)} :=v$.
\item  Draw a random number $s$ from uniform distribution on $[0,1]$. If $s \le  \frac{\pi(v|{\bf m})}{\pi(\textbf{$v^{i}$} | {\bf m})}$ then set $v^{(i+1)}=v$, else set $v^{(i+1)} := v^{(i)}$. 
\item If $i=N$ then stop; else set $i:=i+1$ and go to $2^{nd}$ step.
\end{enumerate}

To improve the efficiency of the Metropolis Hastings algorithm: Delayed Rejection Adaptive Metropolis (DRAM) is proposed \cite {Andrieu,Haario,Rosenthal,Tierney1998}. DRAM is a method that combines two powerful ideas: Delayed rejection (DR) and Adaptive Metropolis (AM).  In DR, if the proposal state is rejected, instead of remaining in the same position a second stage is proposed. The strategy of AM combining to DR method allows the covariance matrix calibrated using {\it n}-samples path of the MCMC chain regardless at which stage these samples are accepted in DR algorithm.

\subsection{Total Variation regularization} \label{Total Variation regularization}
In this tomographic model, we have only a homogeneous object and a sharp boundary that divides the background and the domain of the object. One of the well known methods to solve this problem is total variation (TV) regularization which produces {\it edge-preserving} reconstruction. Let us denote $\Phi = [b_{ij}]$. The solution of TV-regularized is defined by finding the vector $\Phi $ that minimizes the penalty functional:

\begin{equation}\label{eq:TV}
\Phi^{(\alpha)}=\arg\min_{\Phi \in \R^{N^2}} \{\|  \mathcal{A}\Phi-{\bf m}\|_2^2 + \alpha {\text TV}(\Phi)\},
\end{equation}
where $\alpha>0$ is a regularization parameter and TV is a total variation defined as follows.

$$
{\text TV} (\Phi)=\sum_{i=1}^N \sum_{j=2}^N |b_{ij} - b_{i(j-1)}| + \sum_{i=2}^N \sum_{j=1}^N |b_{ij} - b_{(i-1)j}|.
$$

To apply that in the computation, 2D adaptation of the quadratic programming approach is implemented. See \cite{Mueller} Section 6.2.
\subsection{Optimal TV parameter choice} \label{Optimal TV}
In this subsection, we discuss about how to choose an optimal regularization parameter for TV in terms of a piecewise constant image.
Because TV is the baseline method in this study, we want it to perform optimally. Therefore, we choose the regularization parameter in an unrealistically effective way, using the knowledge of the true (simulated) attenuation coefficient.

In the case of simulated data, we implement the thresholding method to yield the result as a binary image by choosing the optimal regularization parameter $\alpha$ and the attenuation threshold $\beta$, where \[ 0<\beta\leq  \max_{ij}   b_{ij}^{{(\alpha)}}. \]

Let us write $\Phi^{(\alpha)} = [b_{ij}^{(\alpha)}]$, where $\alpha$ is the optimal regularization parameter. Define $\Phi^{(\alpha,\beta)} = [b_{ij}^{(\alpha,\beta)}]$ as thresholded TV reconstruction, where

\begin{equation}\label{eq:b}
b_{ij}^{(\alpha,\beta)}=
\begin{cases}
0,    \mbox{ if  } b_{ij}^{(\alpha)}<\beta\\
\beta,   \mbox{ if  } b_{ij}^{(\alpha)}\geq\beta.
\end{cases}
\end{equation}

The error of the reconstructions is computed as follows.
Denote $O$ as the image of the target 2D object and $O^{\mbox{\tiny rec}} $ as the image of the reconstruction.  Set $O\backslash O^{\mbox{\tiny rec}} $ for points that belong to the original object but not to the reconstruction and $O^{\mbox{\tiny rec}} \backslash O $ for points that belong to the reconstruction but not to the original object. The relative error in the reconstruction is written as
\begin{equation}\label{RE}
\frac{(\text{area}(O\backslash O^{\mbox{\tiny rec}} ) + \text{area}(O^{\mbox{\tiny rec}} \backslash O ))}{\text{area }(O)} 100\% . 
\end{equation}

For every $\alpha_i$, the relative error is measured for all $\beta_j$. Eventually, we pick the regularization parameter $\alpha$ that contains the minimum error from each thresholding value $\beta$.

\subsection{Simulated phantoms} \label{Simulated phantoms}
\begin{figure}
\begin{picture}(100,90)
\put(50,0){\includegraphics[width=3.25cm]{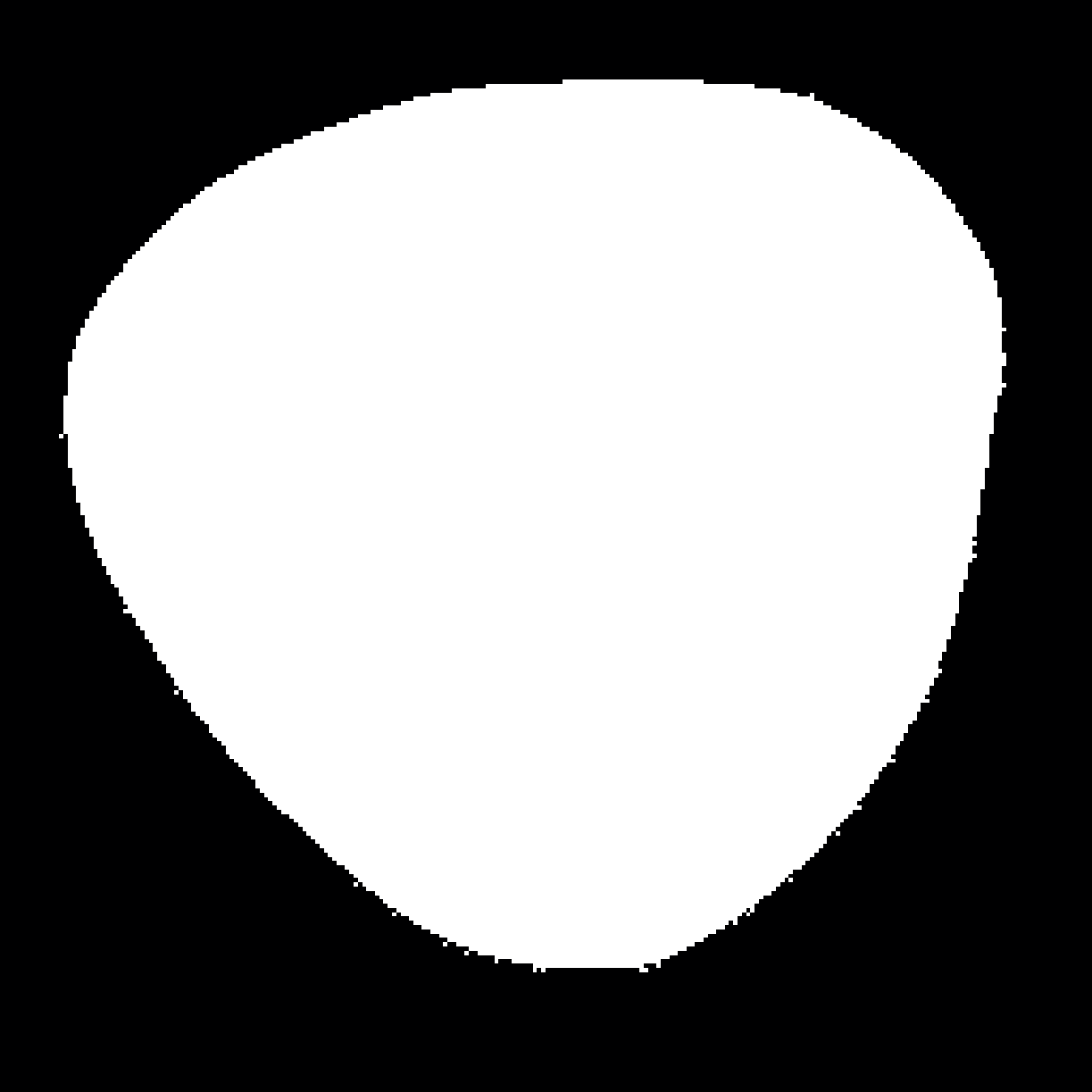}}
\put(200,0){\includegraphics[width=3cm]{original.eps}}
\end{picture}
\caption{\label{fig:simulated phantoms} Left: Simulated phantom $\Omega_1$. Right: Simulated phantom $\Omega_2$.}
\end{figure}

In this part, two simulated phantoms are presented. We build images with one homogeneous convex shape and one homogeneous non-convex shape as in Figure \ref{fig:simulated phantoms}.
 
Each phantom has the image resolution $256 \times 256$ and they are built without using NURBS to avoid an inverse crime.
The images have only two values: 0 for the background and 0.027 for the inner shape.

\subsection{Physical phantoms}\label{Physical phantoms}
\begin{figure}
\begin{picture}(100,90)
\put(50,0){\includegraphics[width=3.25cm]{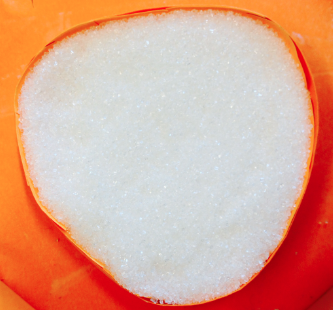}}
\put(200,0){\includegraphics[width=3.2cm]{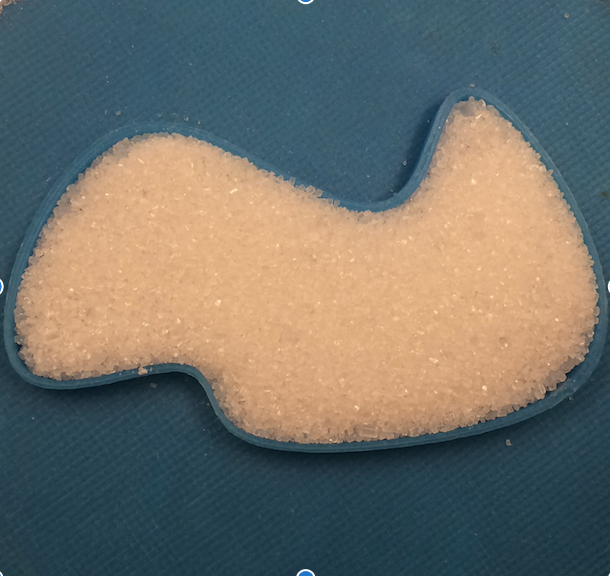}}
\end{picture}
\caption{\label{fig:Physical phantoms} Left: Physical phantom $\Omega_1$ made by cardboard. Right: 3D printed physical phantom $\Omega_2$.}
\end{figure}

In this paper, we also present the reconstruction of real objects with the same shape as the simulated phantoms. One of the physical phantom $\Omega_1$ is simply made by cardboard while another phantom $\Omega_2$ is plastic-printed using 3D printer machine. Both are then filled by white crystal sugar as shown in Figure \ref{fig:Physical phantoms}.

\section{Results} \label{Results}
Throughout the chapter, we use the following measurement setup: sparse-angle (6 angles) from full $360^{\circ}$ angles of {\it fan beam} geometry are presented. 

\subsection{Simulated phantom reconstruction} \label{Simulated phantom reconstruction}

The projection data of both simulated phantoms, $\Omega_1$ and $\Omega_2$, is corrupted by $0.1\%$ Gaussian noise each. In the real projection data, we normalize it as $\mathbf{max}(\mathbf{log}(m))-\mathbf{log}(m)$. Based on the {\it a priori} information that the target is solid homogeneous object, we set the small values (less than $0.1$) in the $m$ to be zero. This intends to be an air background.

Prior to implementing the MCMC algorithm, some parameters need to be fixed. Based on the experiments, to create the perfectly target shapes of $\Omega_1$ and $\Omega_2$ in terms of NURBS curve, the minimum number of control points that are needed are 6 and 12, respectively. Those numbers of control points are then proceeded as our parameters of interest for both simulated and real problems. The periodic uniform knot vector is used and the weights are set equally the same for all control points. As discussed in Subsection \ref {Tomographic Measurement Model}, the coordinates of the control points are defined in the polar coordinates. Since the parameter of the attenuation value needs to be recovered as well, therefore the total number of parameters of interests  for $\Omega_1$ and $\Omega_2$ are 13 and 25, respectively. In DRAM scheme, as a rule of thumb, the length of non-adaptation period in low dimensional problems, $n_0$, is fixed to be $100$.

The simulations are performed with a MATLAB code on the machine equipped with 2.9 GHz Intel Core i5 CPU and 8 GB memory. The image resolution in simulation run is set to have a size of $64 \times 64$ pixels for the sake of computing speed. A unit disc with the radius $ \frac {2}{5}$ over the resolution width and the attenuation value $5$ are fixed as starting points. A Gaussian distribution with centre at $\widetilde {v} =[\widetilde{r_1}, \widetilde{\theta_1}, ... ,\widetilde{r_{\bf n}},\widetilde{\theta_{\bf n}},\widetilde{c} ]^T$ is set for the prior where {\bf n} is the number of control points. For both targets $\Omega_1$ and $\Omega_2$, $\widetilde{r_i}=32$ for $i=1,...,{\bf n}$ and $\widetilde{c}=0.1$. The angles are set to be $\widetilde{\theta}=0^\circ,60^\circ,120^\circ,180^\circ,240^\circ,300^\circ$ for phantom $\Omega_1$ and $\widetilde{\theta}=0^\circ,30^\circ,60^\circ,90^\circ,120^\circ,150^\circ,$ $180^\circ,210^\circ,240^\circ,270^\circ,300^\circ,330^\circ$ for phantom $\Omega_2$. 

The total number of evaluations is 6\,000\,000. The histograms and the chains for simulated data reconstructions show relatively the same behaviour as those of the physical data. Therefore the figures are omitted and the histograms and the MCMC chains of physical data reconstruction are given in the Subsection \ref{Physical phantom reconstruction}.

\begin{figure}
\begin{picture}(200,200)
\put(60,110){\includegraphics[width=3.25cm]{originalO.eps}}
\put(190,110){\includegraphics[width=3.25cm]{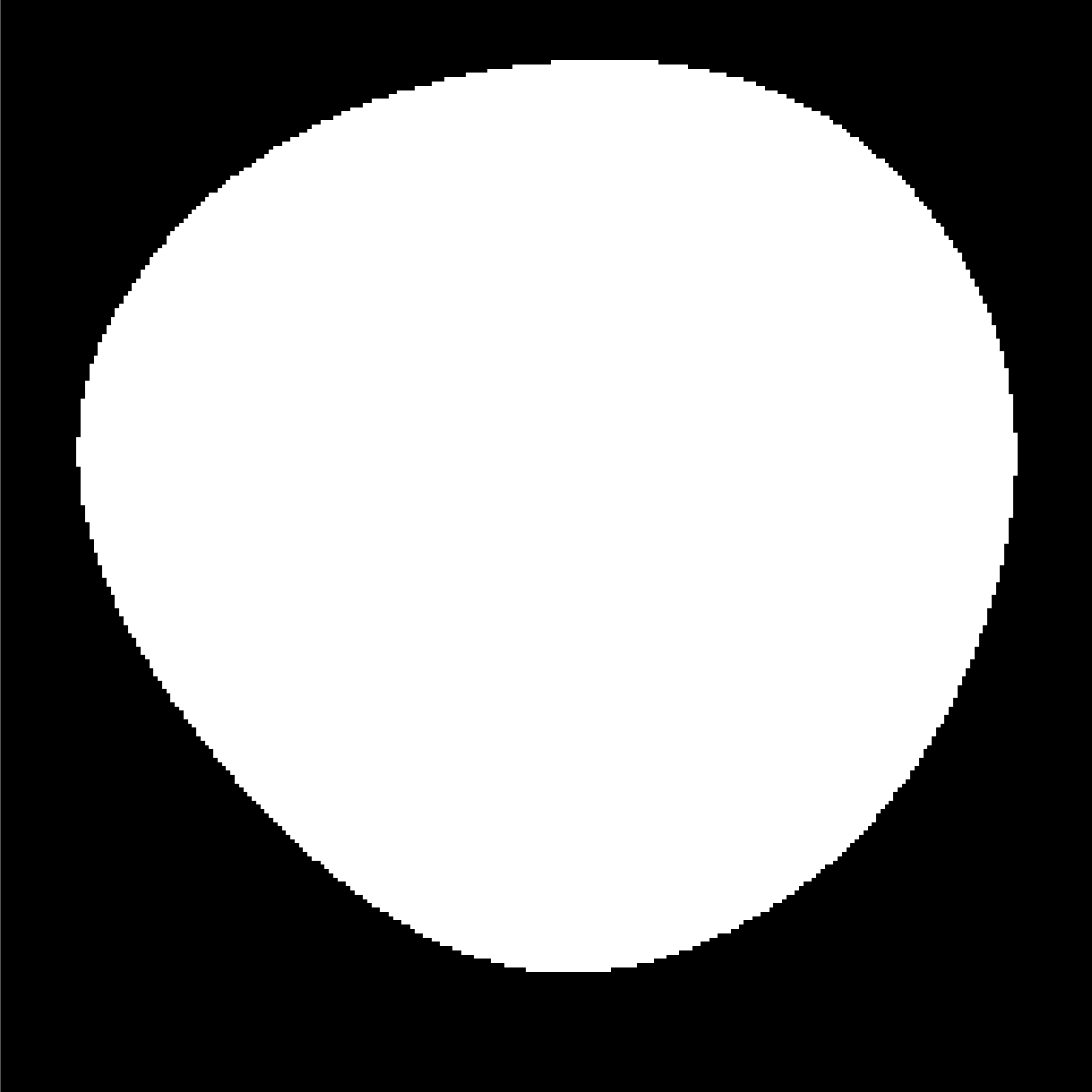}}
\put(60,4){\includegraphics[width=3.25cm]{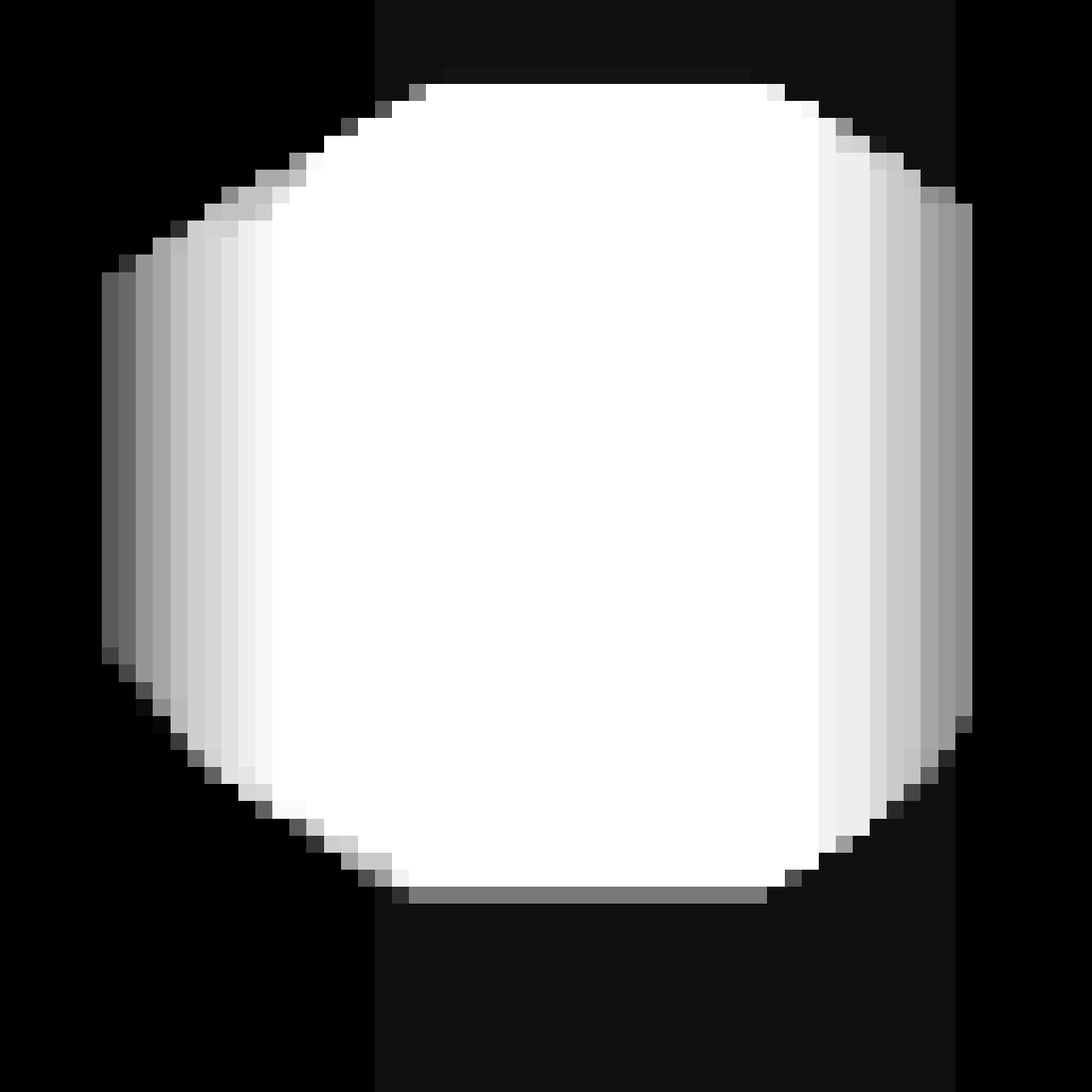}}
\put(190,4){\includegraphics[width=3.25cm]{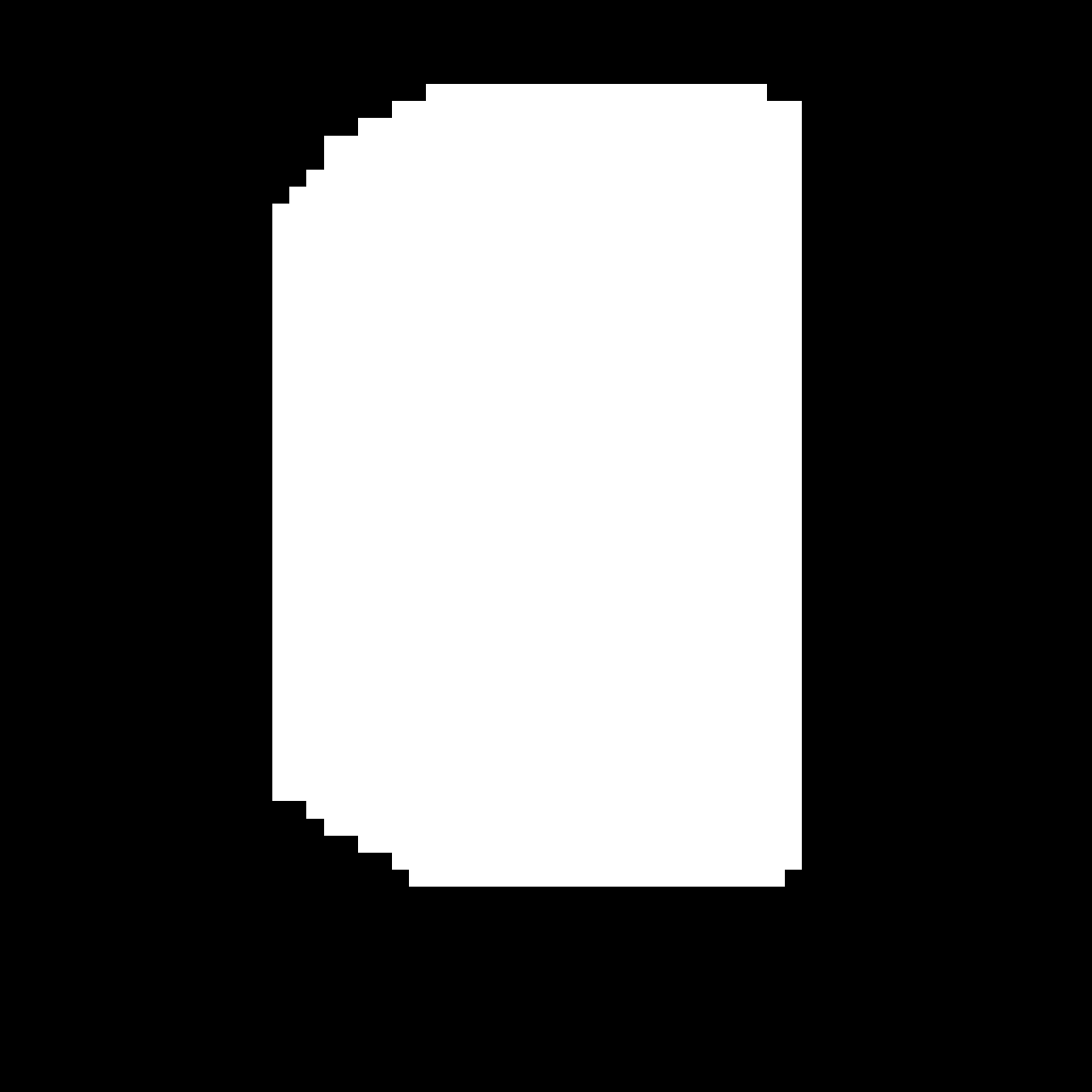}}
\put (40,193){(a)}
\put (170,193){(b)}
\put (40,85){(c)}
\put (170,85){(d)}
\end{picture}
\caption{\label{fig:simulated_1 6 angles} Reconstruction of simulated phantom $\Omega_1$ with 6 projections. (a) The ground truth. (b) NURBS-MCMC reconstruction. (c) TV regularization. See Table~\ref{tab:TV1 data} for the   optimal parameters used. (d) Thresholded TV reconstruction.}
\end{figure}

\begin{figure}
\begin{picture}(200,200)
\put(60,110){\includegraphics[width=3.25cm]{original.eps}}
\put(190,110){\includegraphics[width=3.25cm]{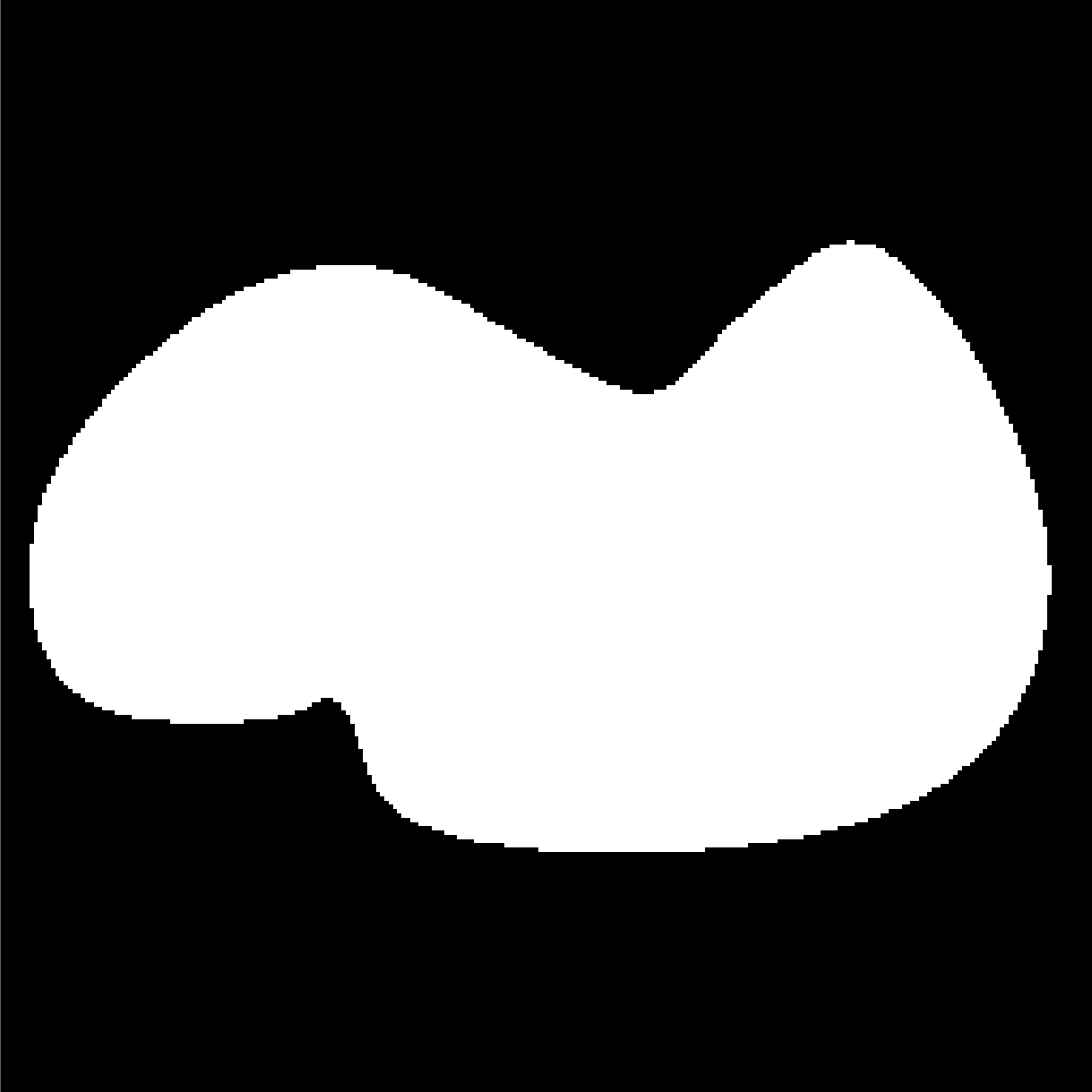}}
\put(60,4){\includegraphics[width=3.25cm]{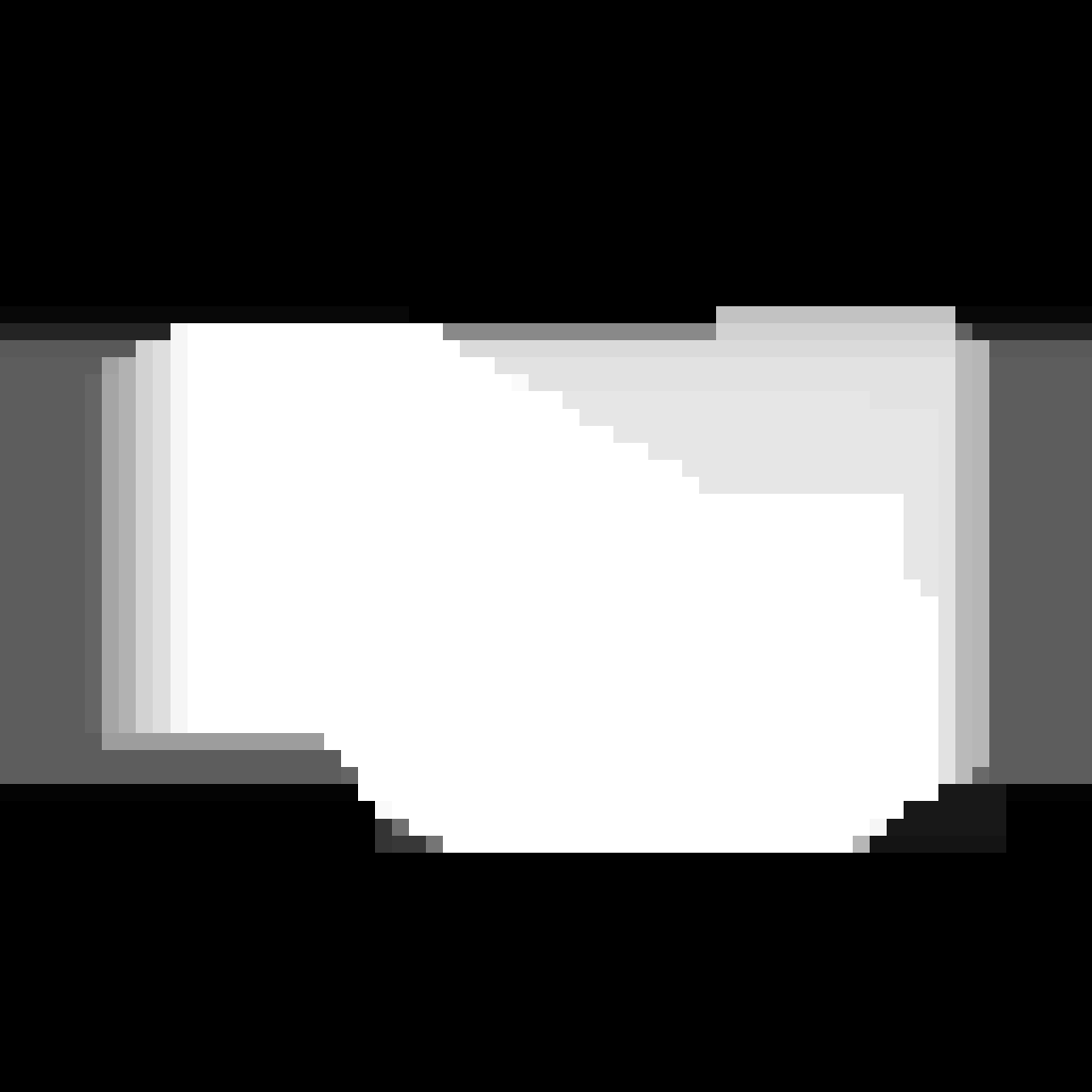}}
\put(190,4){\includegraphics[width=3.25cm]{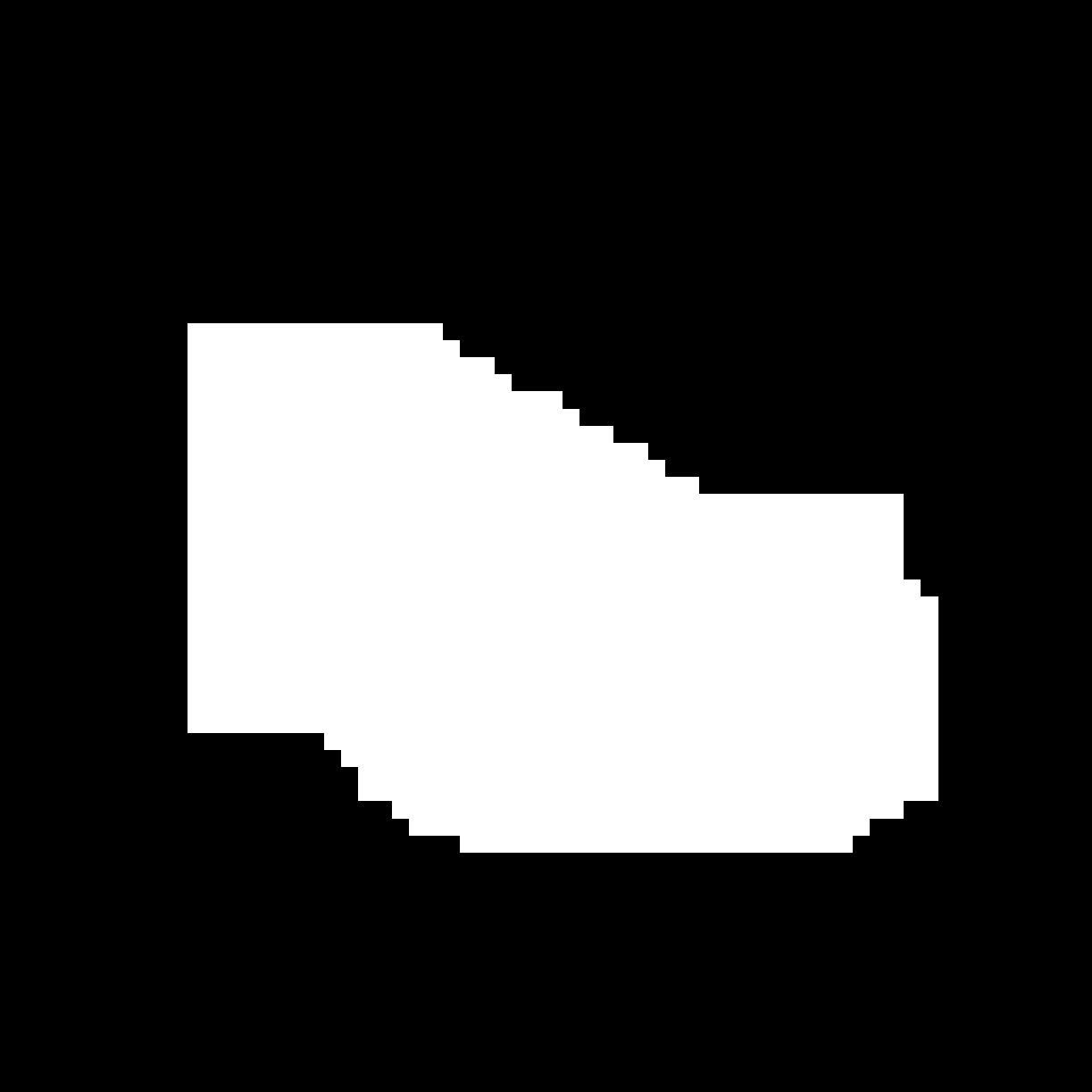}}
\put (40,193){(a)}
\put (170,193){(b)}
\put (40,85){(c)}
\put (170,85){(d)}
\end{picture}
\caption{\label{fig:simulated_2 6 angles} Reconstruction of simulated phantom $\Omega_2$ with 6 projections. (a) The ground truth. (b) NURBS-MCMC reconstruction. (c) TV regularization. See Table~\ref{tab:TV1 data} for the optimal parameters used. (d) Thresholded TV reconstruction.}
\end{figure}

For comparison, the TV regularizations for the 2D tomographic case using quadratic programming as well as the thresholded-TV reconstructions are presented. See Figures~\ref{fig:simulated_1 6 angles} and \ref{fig:simulated_2 6 angles}. The optimal TV parameter choice is calculated as discussed in Subsection \ref{Optimal TV}. In the thresholded-TV reconstructions, we choose the range of regularization parameter range of $10^{-6} \leq \alpha \leq 100$ and the {\it attenuation threshold} value of $0.01 \leq \beta_j \leq 0.03$. The optimal parameters choices for each data are given in Table ~\ref{tab:TV1 data} (Table ~\ref{tab:TV2 data} is the optimal parameters choices for physical problem). The absolute relative errors are presented as well as a comparison to CM estimate of the attenuation value, $c^{\mbox{\tiny CM}}$.

\begin{table}
\caption {The optimal parameters of the best Total Variation reconstructions of simulated phantoms. See Equation \ref{eq:TV} and Equation \ref{eq:b} for the definition of $\alpha$ and $\beta$, respectively.} \label{tab:TV1 data} 
\begin{center} 
\begin{tabular}{| l | c | c | } \hline 
& Optimal $\alpha$&  Optimal  $\beta$   \\ \hline 
Simulated $\Omega_1$&  82.76 & 0.0249  \\ \hline 
Simulated $\Omega_2$ & 82.76  & 0.0225   \\ \hline 
\end{tabular} 
\end{center}
\end{table}

\begin{table}
\caption {The optimal parameters of the best Total Variation reconstructions of physical phantoms. See Equation \ref{eq:TV} and Equation \ref{eq:b} for the definition of $\alpha$ and $\beta$, respectively.} \label{tab:TV2 data} 
\begin{center} 
\begin{tabular}{| l | c | c | } \hline 
& Optimal $\alpha$&  Optimal  $\beta$   \\ \hline 
Physical $\Omega_1$ & 31.03  & 0.0254 \\ \hline 
Physical $\Omega_2$ & 93.10  & 0.0215\\ \hline 
\end{tabular} 
\end{center}
\end{table}

\begin{table}
\caption {Relative errors in reconstructed attenuation values.}\label{tab:attenuation_error} 
\begin{center} 
\begin{tabular}{| l | c | c |} \hline 
 &Thresholded Total Variation & MCMC-NURBS  \\ \hline 
Simulated $\Omega_1$&  7.8\%& 0.37\% \\ \hline 
Simulated $\Omega_2$ & 16.7\% & 0.74\%  \\ \hline 
\end{tabular} 
\end{center}
\end{table}

\begin{table}
\caption {Errors in reconstructed shape of simulated phantoms. The error is defined in Equation \ref{RE}.} \label{tab:shape_error1} 
\begin{center} 
\begin{tabular}{| l | c | c |} \hline 
 &Total Variation & MCMC-NURBS  \\ \hline 
Simulated $\Omega_1$&  32.83\% & 7.48\% \\ \hline 
Simulated $\Omega_2$ & 56.93\%  & 3.9\%  \\ \hline 
\end{tabular} 
\end{center}
\end{table}

\begin{table}
\caption {Errors in reconstructed shape of physical phantoms. The error is defined in Equation \ref{RE}.} \label{tab:shape_error2} 
\begin{center} 
\begin{tabular}{| l | c | c |} \hline 
 &Total Variation & MCMC-NURBS  \\ \hline 
Physical $\Omega_1$ & 50.95\%  & 11.67\% \\ \hline 
Physical $\Omega_2$ & 15.7\%  & 5.95\%\\ \hline 
\end{tabular} 
\end{center}
\end{table}

The simulation time for each run is given in Table~\ref{tab:time}.

\begin{table}
\caption {Computation time (in hours). In NURBS-MCMC, the computation time is for 1\,000\,000 iterations. In Thresholded-TV, the computation time is for evaluating 30 regularization parameters for reconstruction. } \label{tab:time} 
\begin{center} 
\begin{tabular}{| l | c | r |} \hline 
&  $\Omega_1$ & $\Omega_2$ \\ \hline 
NURBS-MCMC &  10  & 12  \\ \hline 
Thresholded-TV & 0.67 & 0.67 \\ \hline 
\end{tabular} 
\end{center}
\end{table}

\clearpage

\subsection{Physical phantom reconstruction} \label{Physical phantom reconstruction}
\begin{figure}[h!]
\begin{picture}(150,220)
\put(30,-10){\includegraphics[width=11cm]{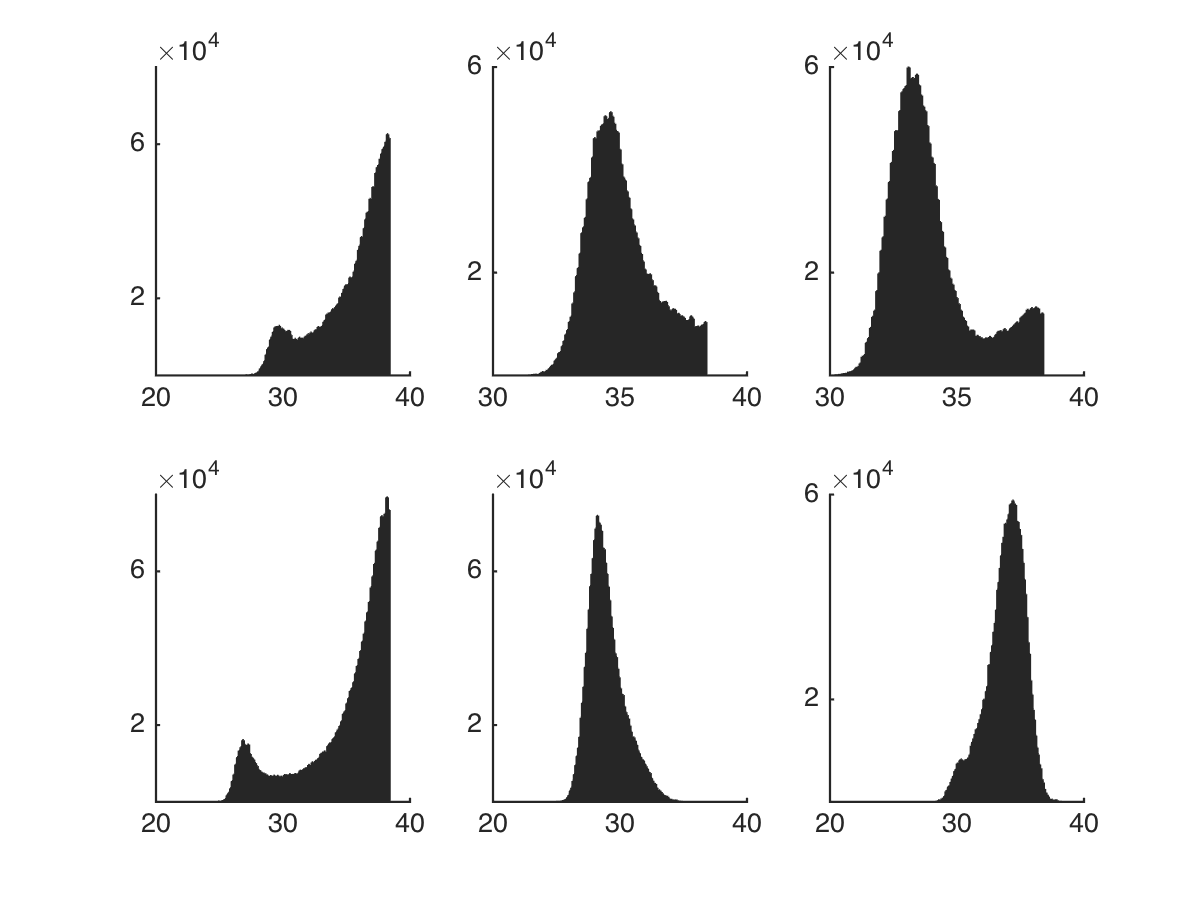}}
\end{picture}
\caption{\label{fig:hist_R_O6r} Histograms of the 1-d marginal posterior distribution of radii in the MCMC chain. This example is related to the physical measurement of $\Omega_1$.}
\end{figure}

\begin{figure}[h!]
\begin{picture}(150,250)
\put(-40,0){\includegraphics[width=17cm]{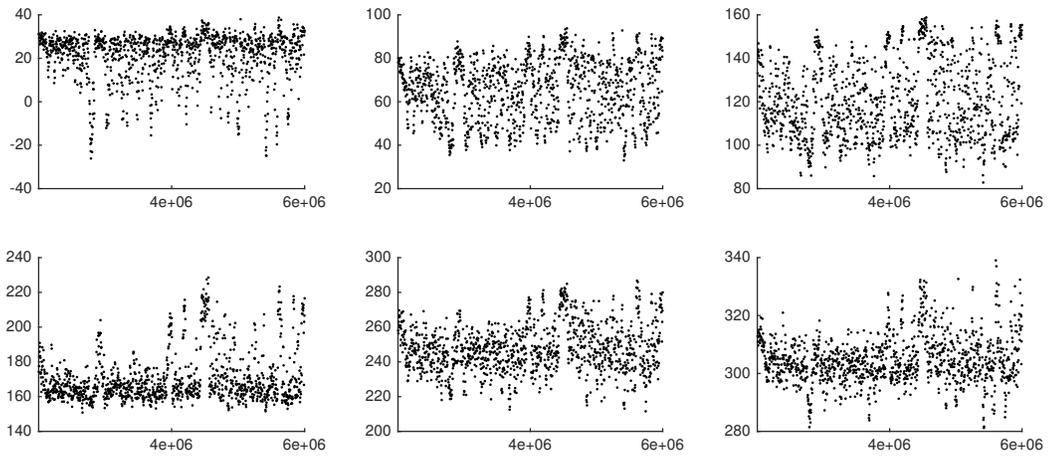}}
\end{picture}
\caption{\label{fig:traceTheta_O6r} The MCMC chains of the angles. This example is related to the physical measurement of $\Omega_1$.}
\end{figure}

\begin{figure}[h!]
\begin{picture}(150,240)
\put(-65,5){\includegraphics[width=19cm]{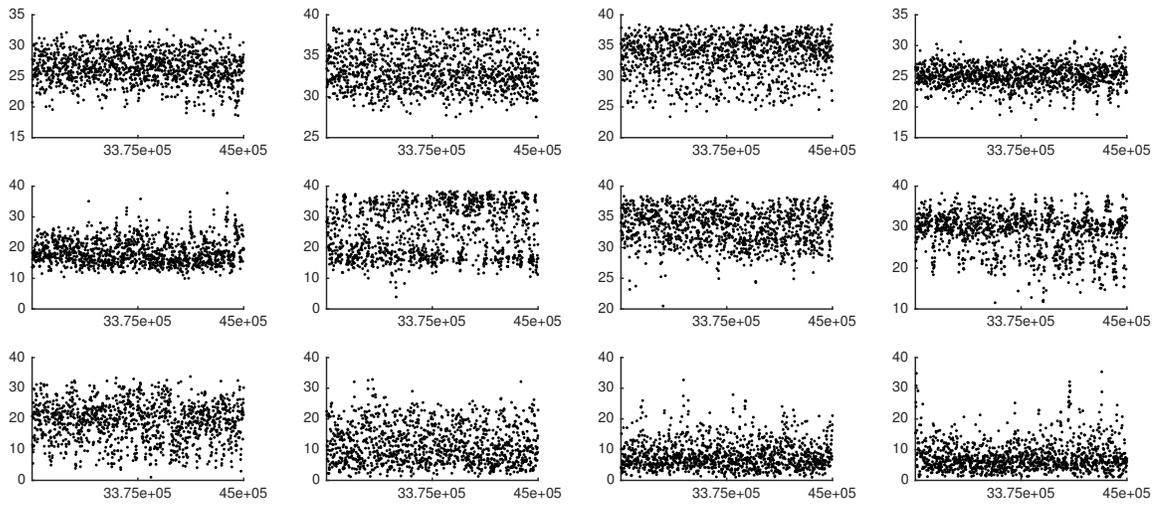}}
\end{picture}
\caption{\label{fig:traceR_B6r} The MCMC chains of the radii. This example is related to the physical measurement of $\Omega_2$.}
\end{figure}

\begin{figure}[h!]
\begin{picture}(150,250)
\put(-10,-5){\includegraphics[width=14cm]{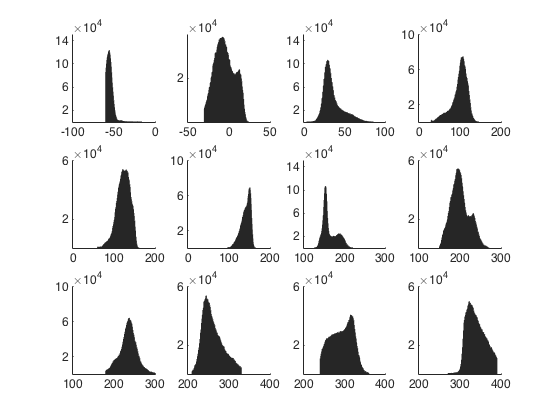}}
\end{picture}
\caption{\label{fig:hist_theta_B6r} Histograms of the 1-d marginal posterior distribution of angles in the MCMC chain. This example is related to the physical measurement of $\Omega_2$.}
\end{figure}

\begin{figure}[h]
\begin{picture}(200,200)
\put(60,110){\includegraphics[width=3.25cm]{originalO.eps}}
\put(190,110){\includegraphics[width=3.26cm]{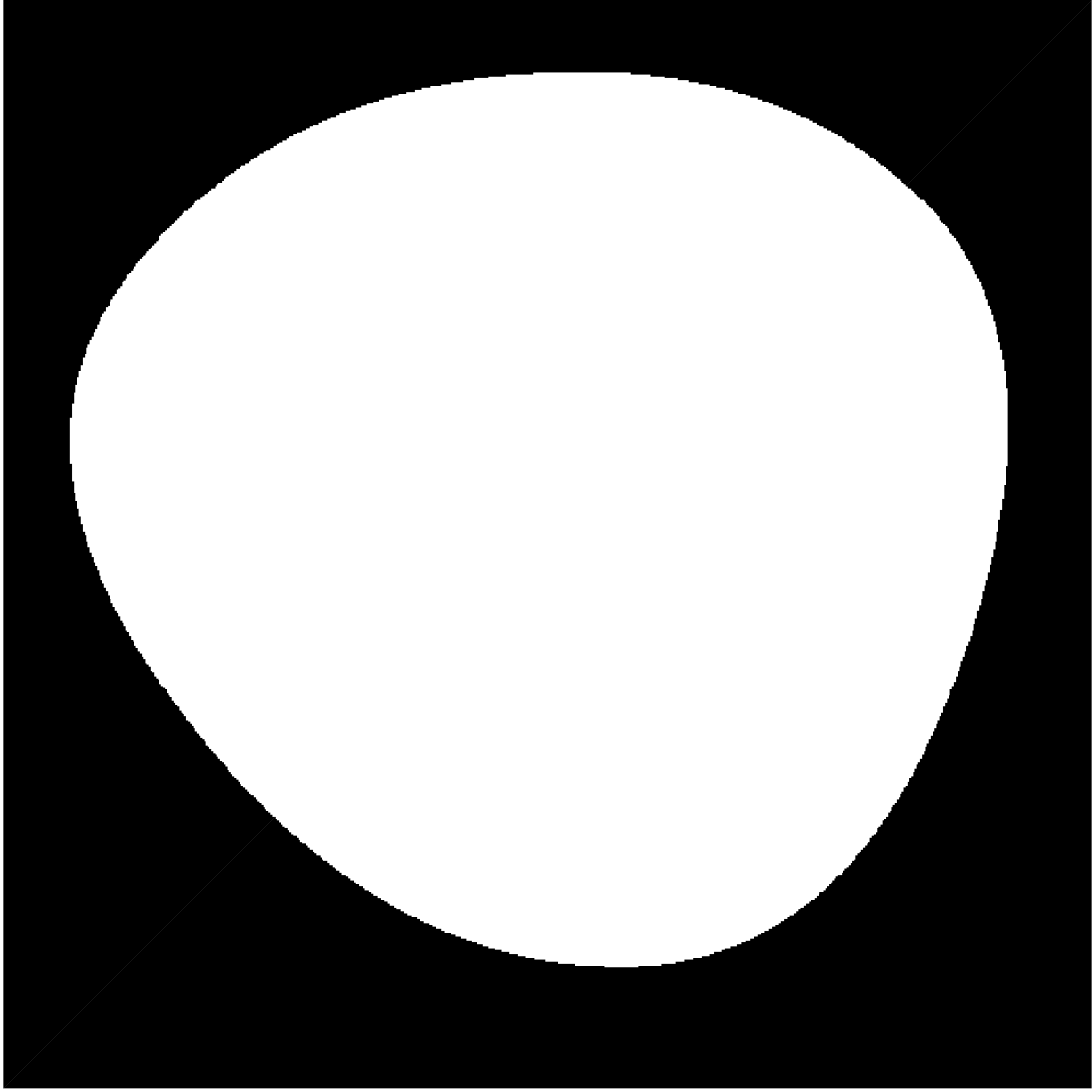}}
\put(60,4){\includegraphics[width=3.25cm]{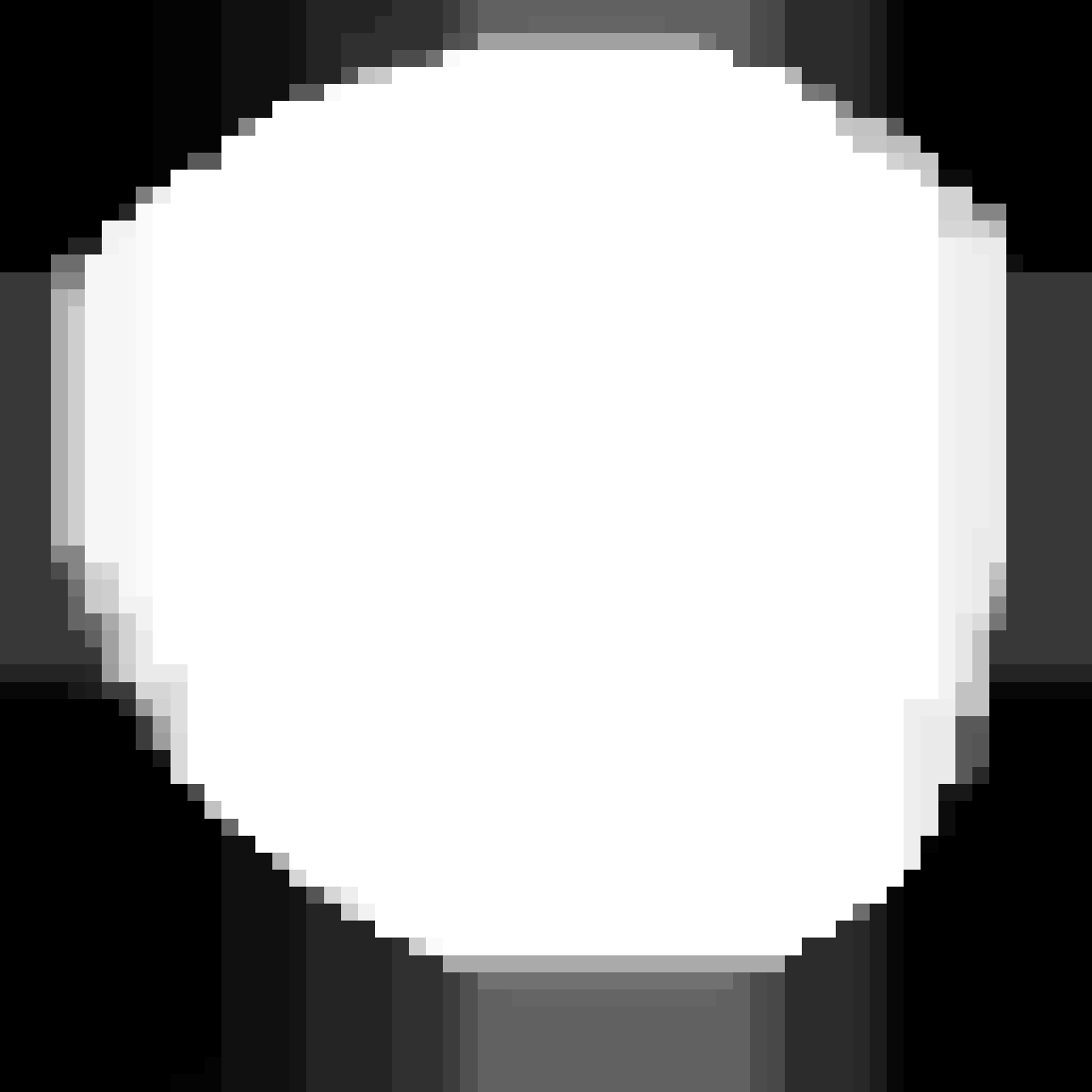}}
\put(190,4){\includegraphics[width=3.25cm]{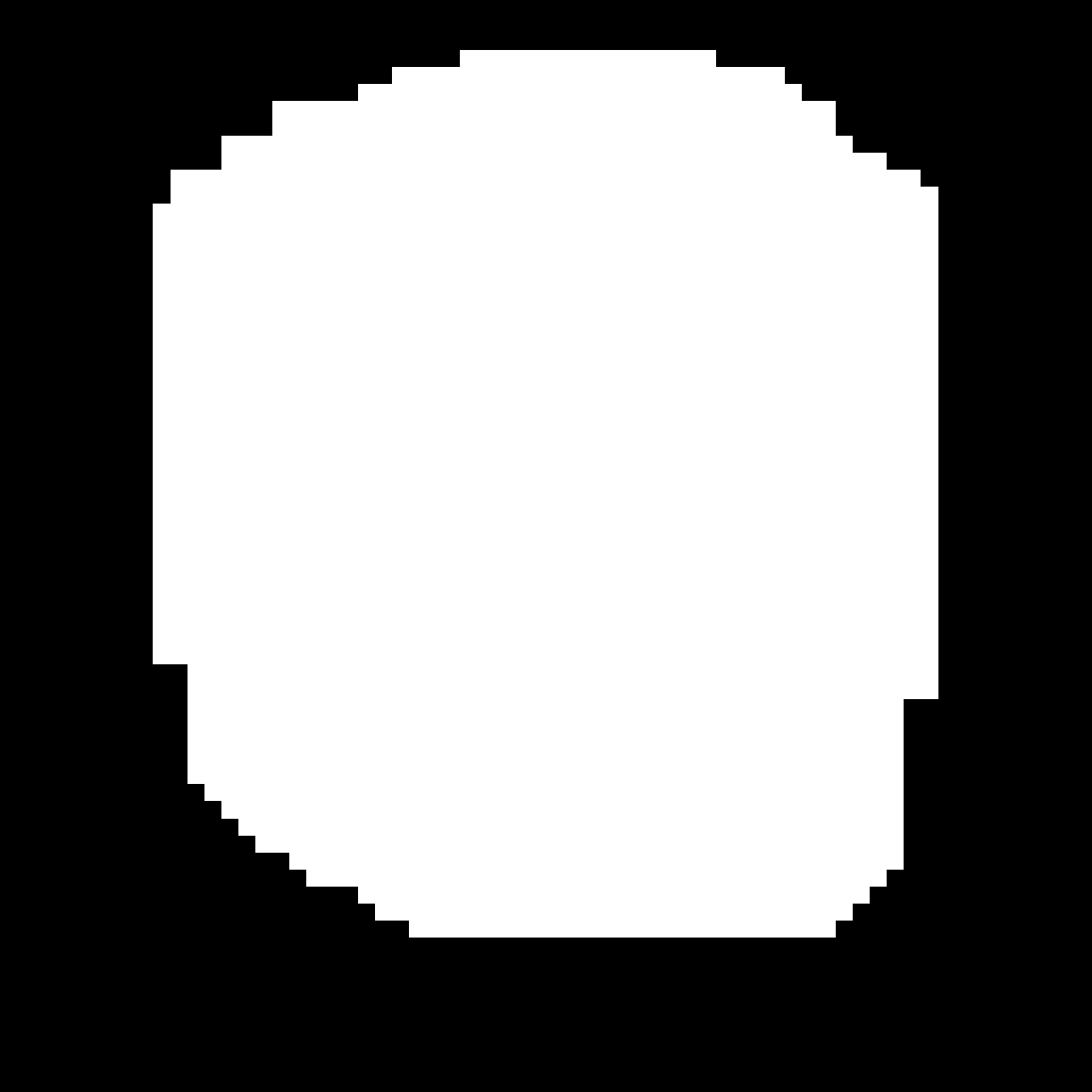}}
\put (40,193){(a)}
\put (170,193){(b)}
\put (40,85){(c)}
\put (170,85){(d)}
\end{picture}
\caption{\label{fig:real_1 6 angles} Reconstruction of physical phantom $\Omega_1$ with 6 projections. (a) The ground truth. (b) NURBS-MCMC reconstruction. (c) TV regularization. See Table~\ref{tab:TV2 data} for the optimal parameters used. (d) Thresholded TV reconstruction.}
\end{figure}

\begin{figure}
\begin{picture}(200,200)
\put(60,110){\includegraphics[width=3.25cm]{original.eps}}
\put(190,110){\includegraphics[width=3.25cm]{ABimg_blue6real.eps}}
\put(60,4){\includegraphics[width=3.25cm]{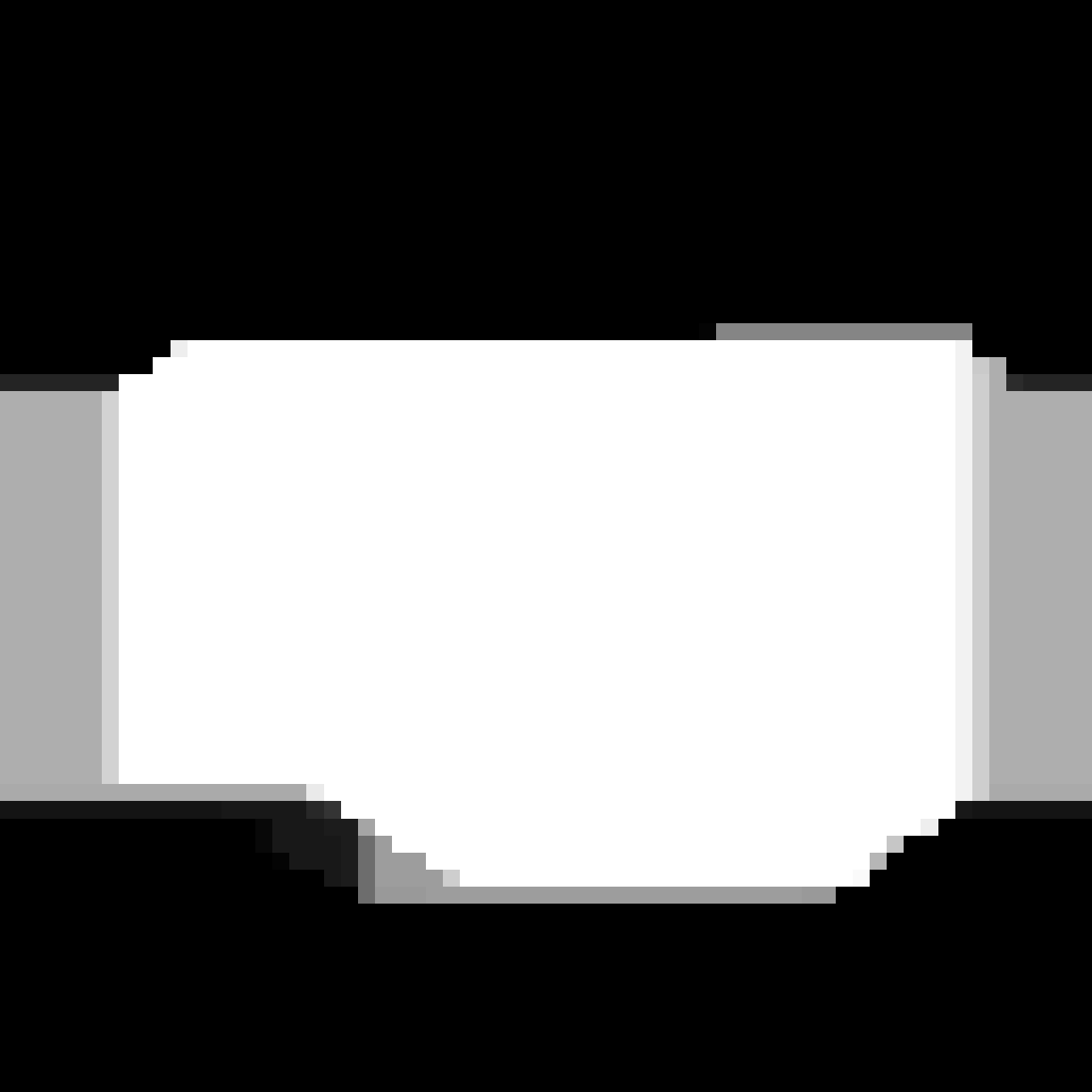}}
\put(190,4){\includegraphics[width=3.25cm]{6Br_last_thresholdedTV_.eps}}
\put (40,193){(a)}
\put (170,193){(b)}
\put (40,85){(c)}
\put (170,85){(d)}
\end{picture}
\caption{\label{fig:real_2 6 angles} Reconstruction of physical phantom $\Omega_2$ with 6 projections. (a) The ground truth. (b) NURBS-MCMC reconstruction. (c) TV regularization. See Table~\ref{tab:TV2 data} for the   optimal parameters used. (d) Thresholded TV reconstruction.}
\end{figure}

The geometry of tomography projection of the real data is given in Subsection \ref{Tomographic Projection Data}. The construction of the prior knowledge is similar with simulated case. The total number of evaluations is 4\,500\,000. Figure \ref{fig:hist_R_O6r} and Figure \ref{fig:hist_theta_B6r} present the histograms of the 1-d marginal posterior distribution of radii from physical data phantom $\Omega_1$ and angles from physical data $\Omega_2$, respectively. The histograms show the distribution of the values of the samples in the MCMC chain. The chains of the angles from physical data phantom $\Omega_1$ and  the chains of the radii from physical data phantom $\Omega_2$ are presented in Figure \ref{fig:traceTheta_O6r} and Figure \ref{fig:traceR_B6r}. All the figures are presented after omitting the {\it burn-in} period. The TV regularizations for the 2D tomographic case using quadratic programming as well as the thresholded-TV reconstructions are presented in Figures~\ref{fig:real_1 6 angles} and \ref{fig:real_2 6 angles}.

\begin{figure}[h!]
\begin{picture}(150,220)
\put(0,-2){\includegraphics[width=6cm]{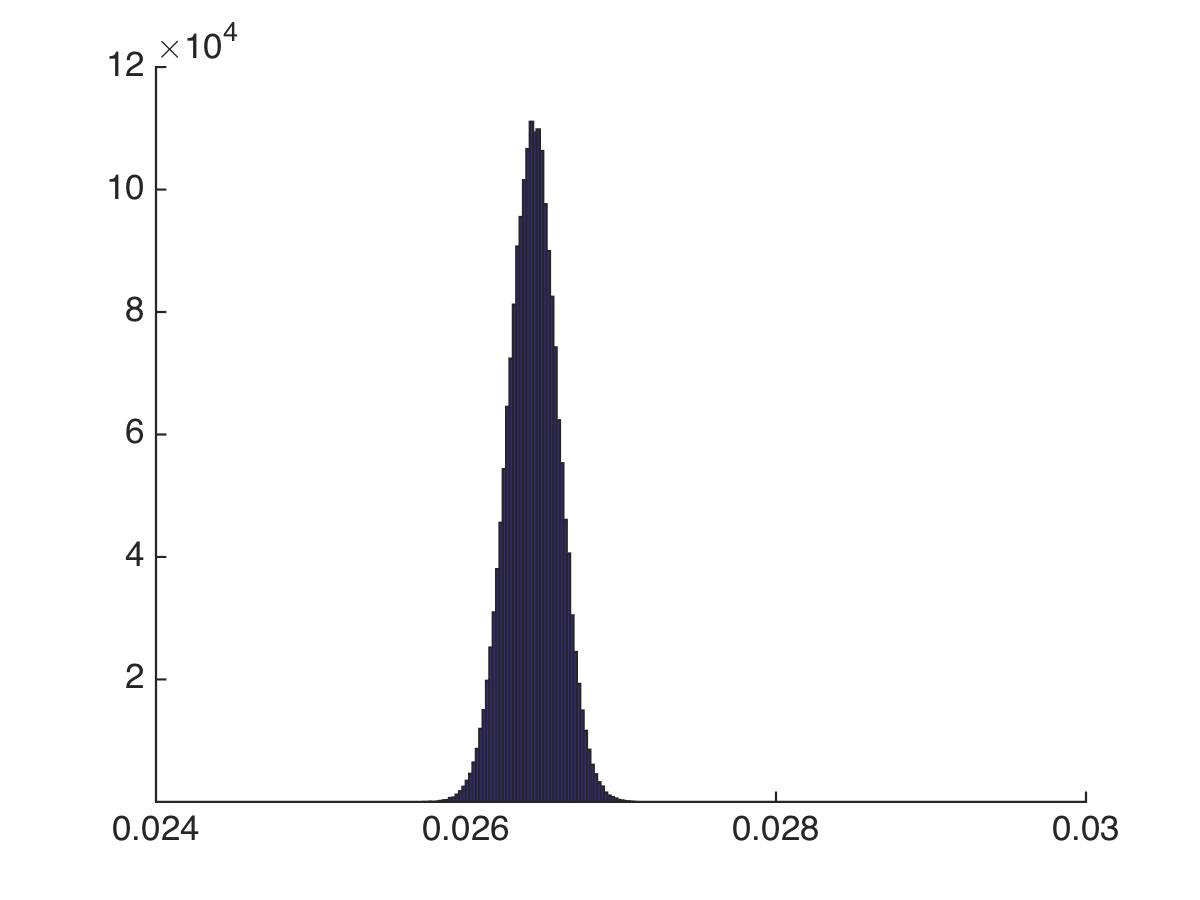}}
\put(170,-2){\includegraphics[width=6cm]{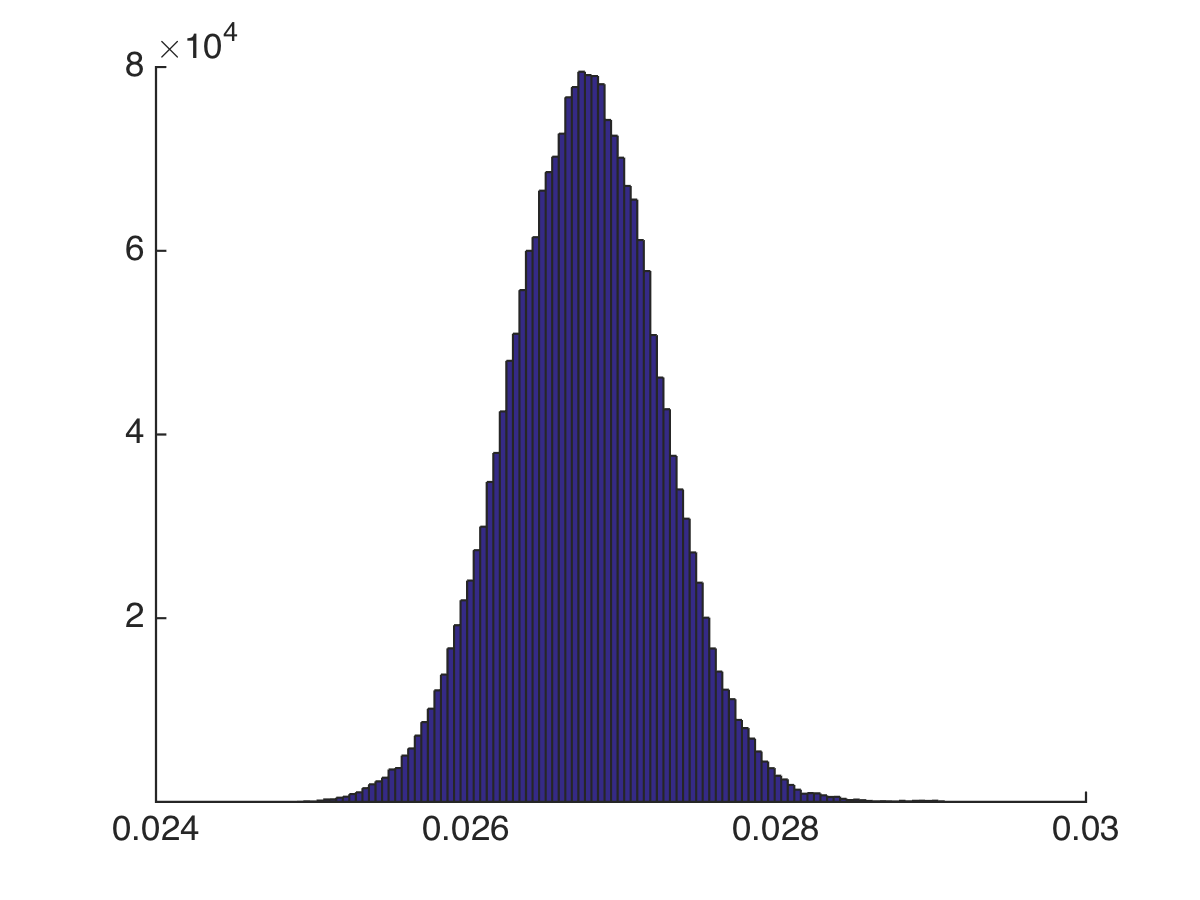}}
\end{picture}
\caption{\label{fig:hist_c} Histograms of the 1-d marginal posterior distribution of the attenuation coefficient in the MCMC chain. This example is related to the physical measurement of $\Omega_1$ (left) and the physical measurement of $\Omega_2$ (right).}
\end{figure}

\begin{figure}
\begin{picture}(200,200)
\put(60,110){\includegraphics[width=3.25cm]{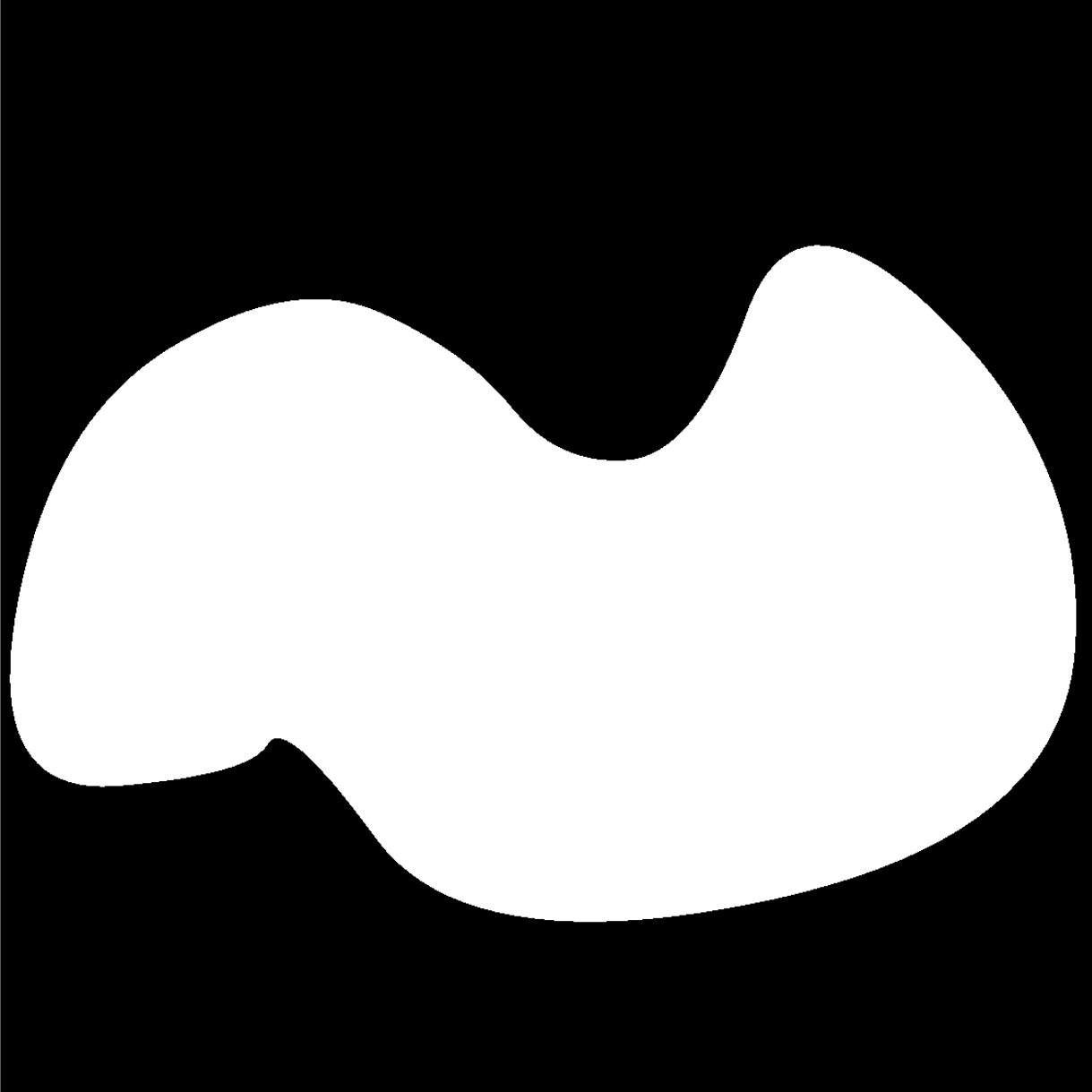}}
\put(190,110){\includegraphics[width=3.25cm]{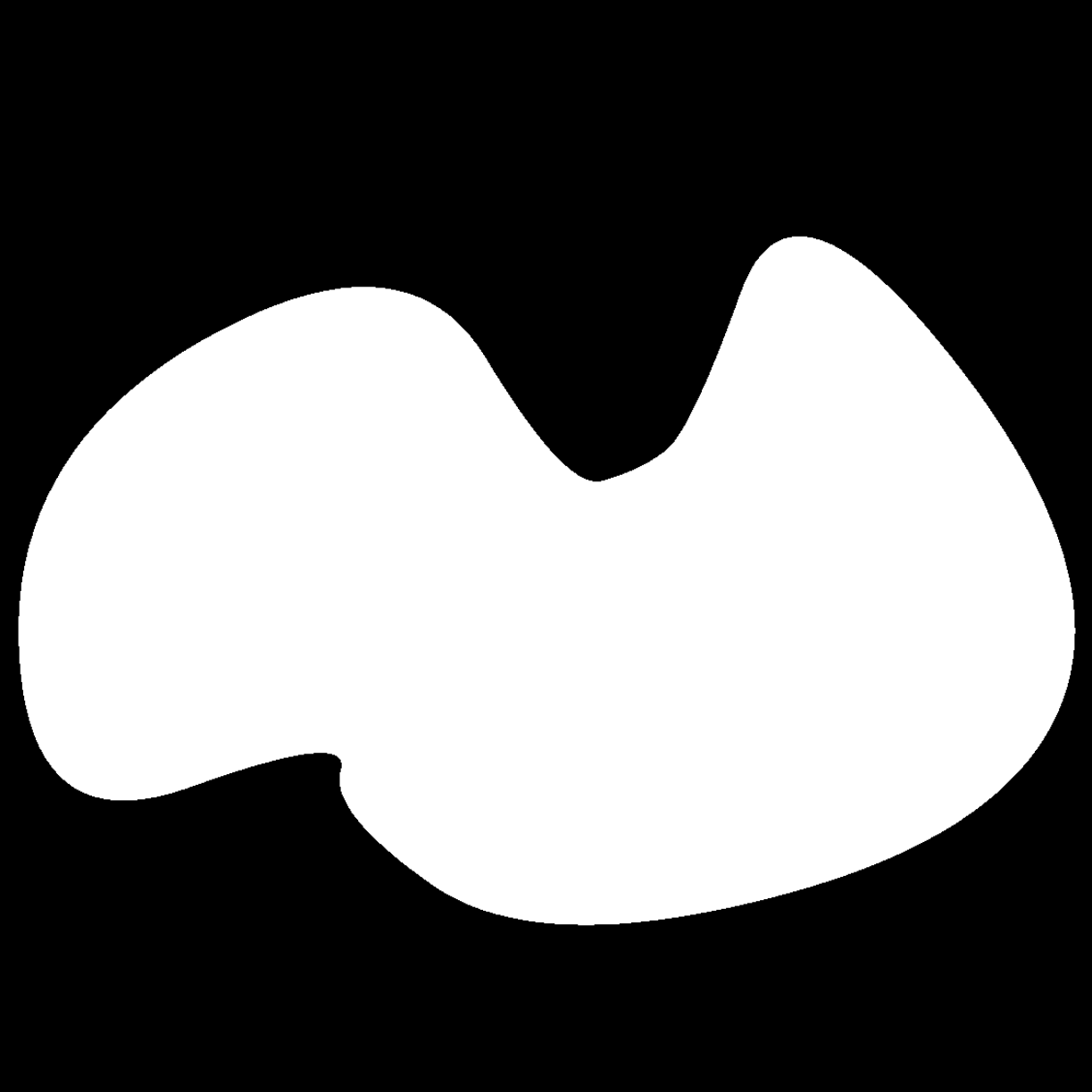}}
\put(60,4){\includegraphics[width=3.25cm]{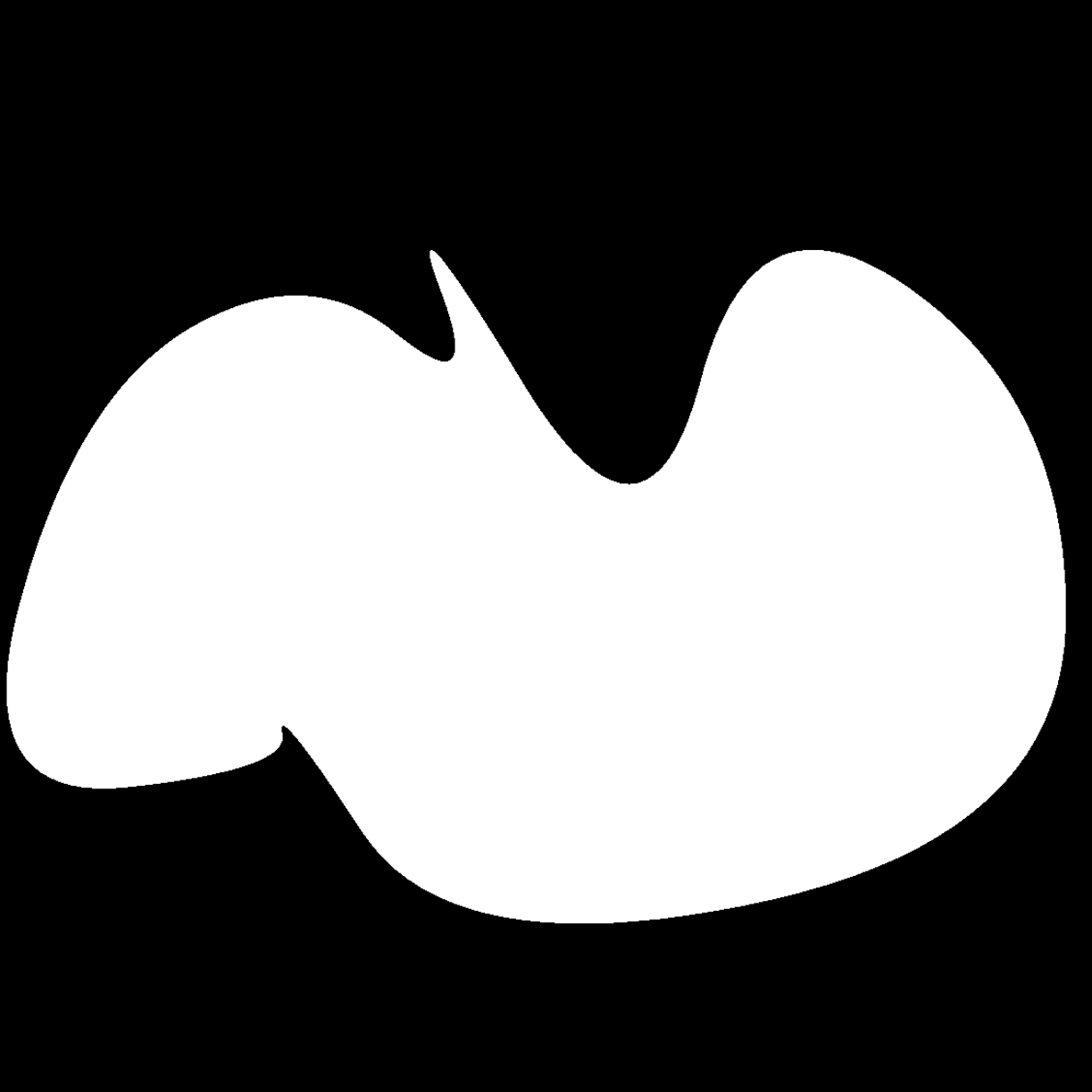}}
\put(190,4){\includegraphics[width=3.25cm]{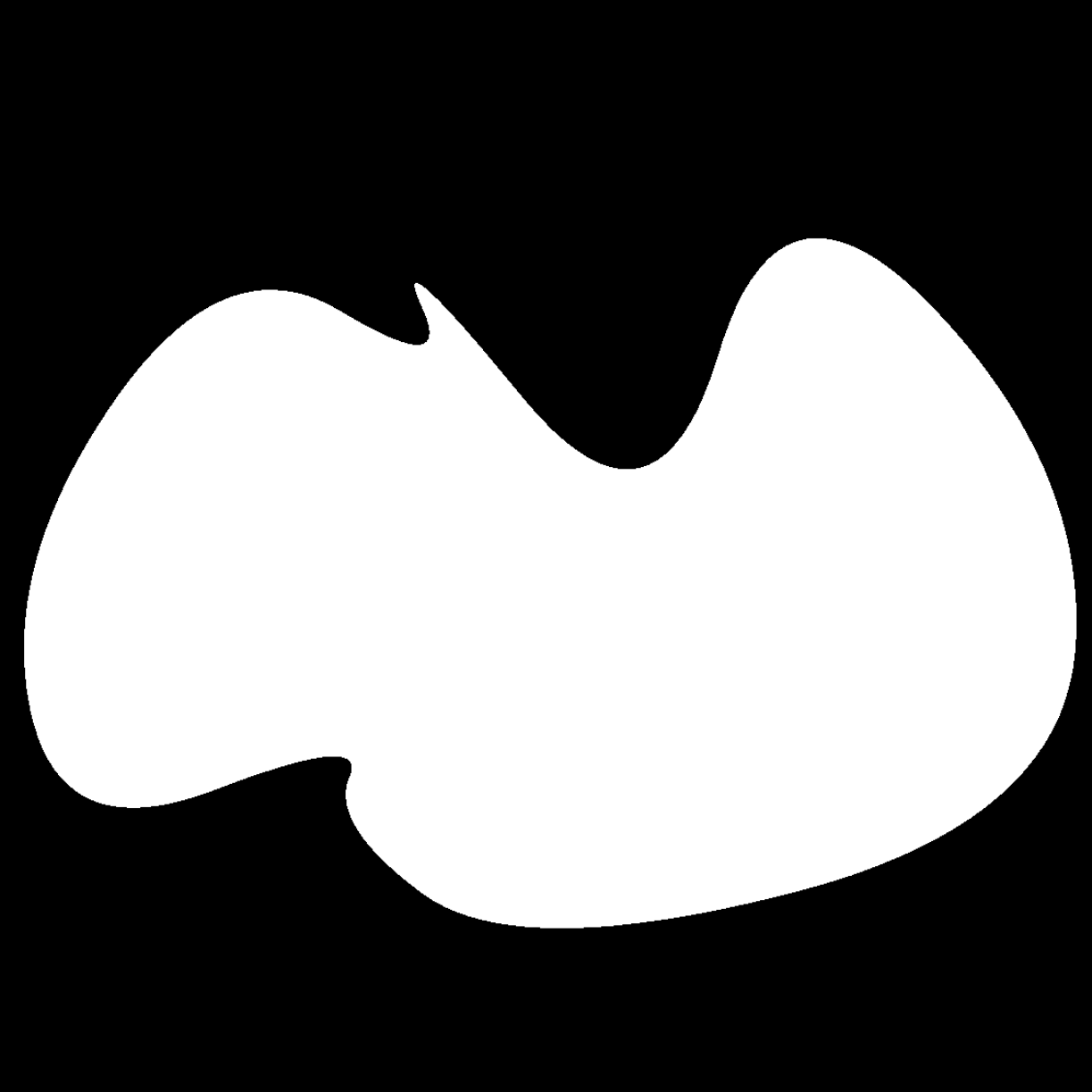}}
\end{picture}
\caption{\label{fig:MAP} Reconstructions of physical phantom $\Omega_2$ with 6 projections which correspond to high posterior values.}
\end{figure}

\clearpage

\section{Discussion}\label{Discussion}

We reconstruct both simulated shapes and outlines of real objects from extremely sparsely collected tomographic data: we only use 6 projection directions (sparsely) spanning the full 360$^\circ$. In our new method we represent the boundaries of homogeneous attenuating objects as a NURBS curve, whose control points are degrees of freedom in the inverse problem. The Bayesian inversion is done by sampling the posterior distribution with an MCMC method. By referring to Section \ref{Tomographic Measurement Model}, another degree of freedom, the attenuation coefficient is set as well as our parameter of interest. The final result, in the form of the NURBS curve, is obtained by computing the conditional mean (CM estimate) of the chains of control points and the attenuation coefficient.

Our baseline method for comparison is Total Variation regularization followed by a thresholding step to yield a homogeneous object. To be as fair as possible to the comparison method, we chose both the regularization parameter and the binary-image threshold optimally.

Let us discuss the features of our Bayesian inversion computation. The MCMC chains of all radii and angles seem to be mixing very well. The convergence is determined by visually looking at one- and two-dimensional chain plots and by Geweke's convergence diagnostic \cite{Brooks}. 

Some of the histograms in Figures \ref{fig:hist_R_O6r} and \ref{fig:hist_theta_B6r} are multimodal and clearly show the non-Gaussian nature of the posterior distribution as a result of our nonlinear NURBS-based parametrization of the unknown. Estimating the full posterior distribution by MCMC allows us to observe nonlinearities such as the multimodality of marginal distributions as shown in Figure \ref{fig:hist_R_O6r}, \ref{fig:hist_theta_B6r} and \ref{fig:hist_c}. We can quantify the uncertainty in the results as well. All of this comes in addition to getting useful reconstructions.

After computing the NURBS curve from the CM estimate of the control points, the result is given as a binary image as an implementation of our tomographic measurement model. Although in the  reconstruction process we used a $64 \times 64$ grid, the final image can be conveniently presented at higher resolution, in this case as $1128 \times 1128$ grid. This is an advantage of the vector-graphic approach we use.

The reconstructions of the simulated and the physical phantoms using the MCMC-NURBS approach produce comparably similar shapes and the attenuation value as the targets. The retrieval of the attenuation parameters of the proposed method for reconstructing the simulated and physical phantoms yields the relative errors of $0.37 \%$ and $1.85\%$, respectively. See Table \ref{tab:attenuation_error}. Curiously, the histograms of attenuation value chains, shown in Figure \ref{fig:hist_c}, are close to Gaussian distributions, in constrast to those of the control points.

As always in Bayesian inversion, the design of the prior distribution is crucial. In this research we spent considerable time for specifying a prior that 
\begin{itemize}
\item is simple to write down as a mathematical formula, and 
\item enables good enough recovery of cavities the non-convex example shapes we use. 
\end{itemize}
The prior we ended up using is a trade-off between several aspects: mathematical simplicity, not too many control points, and flexibility in representing non-convex shapes. 
For example, in Figure \ref{fig:MAP}, the reconstructions from several highest posterior are presented and depicted the change of MCMC chain. In particular, if we look closer to the upper-left of the reconstruction, some of the chain form a spiky shape. This behaviour is allowed by the prior construction which gives the flexibility of the control points to form the cavities. The spiky appearance, of course, can be removed by applying a stricter prior (reduce the possibility of the oscillation more), but consequently, the control points will not form the cavities well. In the case of $\Omega_2$, this is how the cavities are represented: the chains of $r_9$, $r_{10}$, $r_{11}$ and $r_{12}$, combined  with $\theta_9$, $\theta_{10}$, $\theta_{11}$ and $\theta_{12}$, form the bottom cavity. The same holds for $r_3$, $r_{4}$ and $r_{5}$, and for $\theta_3$, $\theta_{4}$ and $\theta_{5}$, in the case of the upper cavity. The prior settings allow the chain to create such cavities. 

A natural follow-up study would be investigating automatical choice of the number of control points. This could be based on a reversible jump MCMC strategy as in \cite{Green}, but we do not discuss such possibilities further here.

In this sparse-data study involving 6 projections only, our proposed method is significantly better than the baseline TV method, as shown in Table \ref{tab:shape_error1} and Table \ref{tab:shape_error2}. Clearly, even our new method still has a room for improvement as the relative errors in shape reconstruction are in the interval 3\%-12\%. However, such improvement must be based on using more {\em a priori} information since the measurement data is very sparse. In this study we do not want to impose more or stricter prior information in the fear of over-fitting our model. In a practical application one might have better and more accurate {\em a priori} information about the target, which could then be used for lowering the reconstruction error further.

There are several interesting avenues for further investigation. Namely, MCMC-NURBS tomography can be used for imaging corrosion inside pipelines; in that case the inner boundary of the pipeline is represented as a NURBS curve. Also, one can extend the new method to targets comprising several components. In that case the MCMC algorithm can have random ``birth'' and ``death'' events varying the number of components. See \cite{Andersen} for such a study involving cracks in electrical conductivity.

\section{Conclusions}\label{Conclusions}

Our NURBS-based nonlinear Bayesian inversion performs very well in the tomographic task of recovering homogeneous 2D objects from extremely sparse projection data.  By using advanced MCMC techniques, we have shown that the MCMC approach is computationally feasible. To tackle the heavy computation, further speed-up strategies are available, such as parallellization, improving the choice of initial value and optimizing the covariances of the sampling strategy. 

The results are conveniently in CAD-compatible vector-graphic format. Quantitative comparison to the baseline method, optimally thresholded TV regularization, is favorable to our method. In the case of recovering the attenuation value, NURBS-MCMC delivers results with relative errors one order of magnitude smaller that the baseline method.

\section*{Acknowledgements}

This work was supported by the Academy of Finland through the Finnish Centre of Excellence in Inverse Problems Research 2012--2017, decision number 250215, and V\"{a}is\"{a}l\"{a} Foundation. We warmly thank Stephan Anzengruber for contributing the integrated code for NURBS.

\end{document}